\documentclass[11pt]{article}

\usepackage{color}
\usepackage{amsmath}
\usepackage{comment}
\usepackage{latexsym,amssymb}
\usepackage{graphicx}
\usepackage{subcaption}
\usepackage{stmaryrd}
\usepackage{fixmath}

\newtheorem{lemma}{Lemma}
\newtheorem{theorem}{Theorem}

\newtheorem{definition}{Definition}

\newtheorem{proposition}{Proposition}

\newtheorem{remark}{Remark}

\def\thelemma{\arabic{section}.\arabic{lemma}}
\def\thetheorem{\arabic{section}.\arabic{theorem}}
\def\thecorollary{\arabic{section}.\arabic{corollary}}
\def\thedefinition{\arabic{section}.\arabic{definition}}
\def\theexample{\arabic{section}.\arabic{example}}
\def\theproposition{\arabic{section}.\arabic{proposition}}

\def\theassumption{\arabic{section}.\arabic{assumption}}

\def\theremark{\arabic{section}.\arabic{remark}}

\newcommand{\manualnames}[1]{
\def\thelemma{#1.\arabic{lemma}}
\def\thetheorem{#1.\arabic{theorem}}
\def\thecorollary{#1.\arabic{corollary}}
\def\thedefinition{#1.\arabic{definition}}
\def\theexample{#1.\arabic{example}}
\def\theproposition{#1.\arabic{proposition}}
\def\theassumption{#1.\arabic{assumption}}
\def\theremark{#1.\arabic{remark}}
}

\newcommand{\noi}{\noindent}

\newcommand{\R}{\mathbb{R}}
\newcommand{\N}{\mathbb{N}}

\newcommand{\la}{\lambda}
\newcommand{\sig}{\sigma}

\newcommand{\eps}{\varepsilon}

\newcommand{\al}{\alpha}
\newcommand{\bet}{\beta}
\newcommand{\io}{\iota}
\newcommand{\gam}{\gamma}

\newcommand{\del}{\delta}

\newcommand{\Gam}{\mathnormal{\Gamma}}

\newcommand{\Sig}{\mathnormal{\Sigma}}
\newcommand{\Ph}{\mathnormal{\Phi}}

\newcommand{\Om}{\mathnormal{\Omega}}

\newcommand{\C}{{\mathbb C}}
\newcommand{\D}{{\mathbb D}}
\newcommand{\M}{{\mathbb M}}

\newcommand{\Z}{{\mathbb Z}}

\newcommand{\calF}{{\cal F}}

\newcommand{\calS}{{\cal S}}

\newcommand{\skp}{\vspace{\baselineskip}}

\newcommand{\diag}{{\rm diag}}

\newcommand{\w}{\wedge}
\newcommand{\lt}{\left}
\newcommand{\rt}{\right}

\newcommand{\To}{\Rightarrow}

\newcommand{\iy}{\infty}
\newcommand{\ds}{\displaystyle}

\newcommand{\IA}{\text{\it IA}}

\newcommand{\qed}{\hfill $\Box$}

\begin{document}

\title{Customer-server population dynamics in heavy traffic}

\author{Rami Atar\thanks{Viterbi Faculty of Electrical Engineering,
Technion --- Israel Institute of Technology,
Haifa 32000, Israel.}
\and
Prasenjit Karmakar${}^*$
\and 
David Lipshutz${}^*$\thanks{Currently at the Center for Computational Biology, Flatiron Institute, New York City, NY 10010, USA.}
}

\date\today

\maketitle

\begin{abstract}
We study a many-server queueing model with server vacations,
where the population size dynamics of servers and customers are coupled:
a server may leave for vacation only when no customers await,
and the capacity available to customers is directly affected
by the number of servers on vacation.
We focus on scaling regimes in which server dynamics and queue dynamics
fluctuate at matching time scales, so that their limiting dynamics are coupled.
Specifically, we argue that interesting coupled dynamics occur in (a) the {\it Halfin-Whitt} regime, (b) the {\it nondegenerate slowdown} regime,
and (c) the intermediate, {\it near Halfin-Whitt} regime; whereas the dynamics asymptotically decouple in the other heavy traffic regimes.
We characterize the limiting dynamics, which are different for each scaling regime.
We consider relevant respective performance measures for regimes (a) and (b) --- namely, the probability of wait and the slowdown.
While closed form formulas for these performance measures have been derived for models that do not accommodate server vacations, it is difficult to obtain closed form formulas for these performance measures in the setting with server vacations.
Instead, we propose formulas that approximate these performance measures, and depend on the steady-state mean number of available servers and previously derived formulas for models without server vacations.
We test the accuracy of these formulas numerically.
\end{abstract}

\section{Introduction}\label{sec1}

Scaling limits for stochastic processing networks with a growing
number of servers is an active research area.
Since the pioneering work of Halfin and Whitt \cite{halfin1981heavy}, the novel scaling that they
introduced has attracted substantial interest, and also inspired the study of various other scaling regimes
in which the number of servers grows to infinity in the limit.
These studies include scaling limit results at the law of large numbers (or fluid) scale,
as well as several distinct frameworks at the central limit theorem (or diffusion) scale.
For a sample of such work, see
\cite{atar2012diffusion},
\cite{ata-gur},
\cite{AMR},
\cite{gam13},
\cite{gam2012},
\cite{gar},
\cite{gup-wal},
\cite{har2004},
\cite{KR1},
\cite{lu2016g},
\cite{pang2009heavy},
\cite{puh-rei},
\cite{van2018},
\cite{whi-04}.

Queueing models with server vacations arise in computer communication systems
and production
engineering, and their mathematical analysis has a long history in the operations research literature; see the earlier survey \cite{doshi1986queueing}, the recent book chapter
\cite[Ch.\ 10]{has16},
and the recent work \cite{lu2016g}, as well as the references therein.
Server vacations occur in models that accommodate primary and secondary classes of customers,
where, from the viewpoint of primary customers, a server working non-preemptively
on secondary customers may equivalently be regarded as if it takes a vacation, as it is not available
during that time.
They also arise for a variety of other reasons, including machine breakdowns and maintenance.
Kella and Whitt \cite{kella1990diffusion} make the distinction between
models in which servers leave for vacations according to the state of the queue and ones
in which vacations are triggered exogenously (see \cite{has16} for various other important distinctions
and classifications of vacation models).
In this paper we consider a model of the former type, where specifically, as in the case
considered in a single-server setting in \cite{kella1990diffusion},
a server may leave only when the queue is empty.
In this case there is an interplay between the population size dynamics of servers and that of customers:
the vacations are triggered by the state of the queue, and the
queue dynamics is affected by the number of available servers.
We study these dynamics at the diffusion scale in a class of heavy traffic many-server
regimes, focusing on
regimes in which the population size dynamics
of customers and of servers fluctuate on the same time scale,
so that the equations describing limiting dynamics remain coupled.
Our first main contribution is to show that diffusion limits can indeed capture
such coupled dynamics, and to classify regimes where it occurs.

To put these results in context
some background on classification of heavy traffic regimes is necessary.
A formulation of a continuum of heavy traffic regimes was introduced in \cite{atar2012diffusion},
which contains as special cases the well-known {\it conventional} and {\it Halfin-Whitt} (HW) regimes.
To introduce it, consider the $N$-server queue, let $\al\in[0,1]$ be a given parameter and let $n$ denote a scaling parameter.
Assume that the arrival rate is proportional to $n$ and that
the number of servers $N_n$ is proportional to $n^\al$.
Impose a critical load condition by letting
the total processing rate be nearly equal to the arrival rate.
Then the individual service rate must be proportional to $n^{1-\al}$.
For any $\al$, a diffusion scaled process is obtained by scaling down the queue size by $n^{1/2}$.
In this spectrum of heavy traffic regimes, the two
endpoints, $\al=0$ and $\al=1$ give the conventional regime (with a fixed number of servers)
and, respectively, the HW regime, with $O(n)$ servers and no acceleration
of service times.
A well-known property of the latter regime, that is unique among all heavy traffic regimes,
is that the steady state probability that an arriving customer
waits in the queue is asymptotic to a number strictly between $0$ and $1$.
At the midpoint, $\al=1/2$, one obtains the {\it nondegenerate slowdown} (NDS)
regime. A unique property of it is that the time in queue and the time in service
for a typical customer are of the same order of magnitude.
The slowdown, defined as the ratio of sojourn time and service time
for a typical customer is therefore nondegenerate in this regime (that is, it is asymptotic
to a number strictly between 1 and $\iy$).
The regimes where $\al\in(0,1/2)$ and $\al\in(1/2,1)$ were not given any names
so far in the literature. In this paper we shall refer to them as the {\it near-conventional}
and the {\it near-HW} regimes, respectively.

The aforementioned coupling between the dynamics of server population size and
queue size is argued in this paper to occur for $\al\in[1/2,1]$ but not for $\al\in[0,1/2)$.
That is, it occurs in the HW, the near-HW and the NDS regimes.
The first setting we consider is of exponential vacation lengths. In this case we show
that the scaling limits of the joint dynamics are governed by a pair of coupled one-dimensional equations. In each of the relevant regimes, the set of equations takes a different form:
\begin{itemize}
	\item [(i)] In the HW regime, the scaling limit is governed by a coupled stochastic differential equation (SDE) and ordinary differential equation (ODE) system.
	\item [(ii)] In the near HW regime, the scaling limit is governed by a coupled SDE with reflection (SDER) at zero and ODE system, where the boundary term that constrains the SDER to remain non-negative also appears
	as a term in the ODE.
	\item [(iii)] In the NDS regime, the scaling limit is governed by a coupled SDER and birth-death process (BDP), where the positive jumps are
	driven by the boundary term for the SDER.
\end{itemize}
We also show that for $\al<1/2$, scaling limits do not exhibit coupled dynamics.
Furthermore, a more involved model is treated,
in which vacation lengths follow a phase type distribution.

The second goal of this paper
is to study the effect of server vacations on natural performance measures,
in the special cases of HW and NDS. In these two cases there are performance measures
that are particularly interesting to study.
In the HW regime, it is the {\it probability of wait}, that is, the steady state
probability that an arriving customer has to wait for service. This probability was shown
in \cite{halfin1981heavy}
to converge to a number strictly between $0$ and $1$, and an explicit formula was given
for the limit.
It is not a meaningful performance measure in
any other regime $\al\in[0,1)$, as in these regimes the limit is always $1$.
We are interested in the asymptotics of the probability of wait in presence of server vacations.
The diffusion limit developed here can in principle make it possible to achieve this goal,
however explicit expressions are hard to obtain for the two-dimensional dynamics (i).
Instead, we propose a further approximation based on heuristics.
This gives rise to a formula that is a variant of
the original formula of \cite{halfin1981heavy}.

Similarly, in the case of NDS, a property that distinguishes this regime
from all regimes with $\al\in[0,1/2)\cup(1/2,1]$ is that the {\it slowdown},
defined as the ratio between expected sojourn time and expected service time
in steady state, is asymptotic
to a random variable strictly between $1$ and $\iy$.
A formula for the slowdown asymptotics was provided in \cite{atar2012diffusion}
for the model without vacations,
and it is of interest to explore how it varies in presence of vacations.
Once again, an explicit expression is hard to obtain as it involves
two-dimensional dynamics, and we turn instead to a heuristic argument.
The heuristic gives rise to a variation of the formula from \cite{atar2012diffusion}
obtained by introducing a correction term.

In both cases, we test the proposed formulas numerically.
We provide arguments suggesting that the heuristic formulas are nearly accurate
in specific parameter settings, and these arguments are validated by our numerical
tests. The overall level of accuracy of the heuristic formulas is also discussed.

There have been relatively few results on diffusion limits for queueing systems with
server vacations, especially in the many-server regime.
In terms of dependence of vacations on the state of the queue,
the closest work to ours is the aforementioned \cite{kella1990diffusion}, which addresses a single
server setting. In their model, the server vacations each time the queue becomes empty.
They also study the case that the server vacations according to an exogenous Poisson process
(in which case there may be service interruptions).
In both cases, the heavy traffic scaling limit is a L\'evy processes with a secondary jump input.
In follow up work \cite{kella1991queues} they prove decompositions for the stationary distributions of
these processes.
In the many-server setting, Pang and Whitt \cite{pang2009heavy} prove diffusion limits
in the HW regime in the case of exogenous server interruptions (see also \cite{lu2016g}
for a related work) that simultaneously affect a proportion of the servers.
This is relevant for models in which exogenous events result in a large number of servers being out of service (e.g., system-wide computer crashes).
The primary difference between our work and the works \cite{lu2016g,pang2009heavy}
is that our focus is on server vacations triggered by the state of the queue,
which leads to a server-customer population dynamics.

The organization of this paper is as follows. Below, some mathematical notation used in this paper
is introduced. In \S \ref{sec2}, the main model and scaling regimes are introduced,
and the first main result is stated. Its proof is provided next in \S \ref{sec3}.
An extension of the result to phase type service time distribution is presented and proved
in \S \ref{sec:multi}. In \S \ref{sec:heu}, heuristic formulas are developed for performance
measures in the HW and NDS regimes. Finally, numerical tests of the level of accuracy of these
heuristics are then provided and discussed in the same section.

\subsection*{Notation}

Let $\N=\{1,2,\dots\}$ denote the positive integers and $\N_0=\N\cup\{0\}$.
For $d\in\N$ let $\R^d$ denote $d$-dimensional Euclidean space.
When $d=1$ we suppress the superscript $d$ and write $\R$ for the real numbers.
Let $\R_+=[0,\infty)$ denote the non-negative axis.
For $a,b\in\R$, the maximum [resp., minimum] is denoted by $a\vee b$ [resp., $a\w b$].
For $a\in\R$, the positive [resp., negative] part is denoted by $a^+=a\vee0$ [resp., $a^-=(-a)\vee0$].
For $a\in\R_+$, let $\lceil a\rceil=\min\{n\in\N_0:n\ge a\}$.
For $d\in\N$ and $x,y\in\R^d$, let $x\cdot y$ and $\|x\|$ denote the usual scalar product and $\ell_2$ norm, respectively.
Given a sequence $\{x_n\}$ in $\R_+$ and $\alpha\ge0$, we say $x_n\sim n^\alpha$ if $n^{-\alpha}x_n\to C$ as $n\to\infty$ for some $C\in\R_+$.

For $f:\R_+\to\R^d$, let $\|f\|_T=\sup_{t\in[0,T]}\|f(t)\|$, and, for $\theta\in(0,T)$, let
\[
w_T(f,\theta)=\sup_{0\le s<u\le s+\theta\le T}\|f(u)-f(s)\|.
\]
For a Polish space $\calS$, let $\C_\calS([0,T])$ and $\D_\calS([0,T])$ denote the set of continuous and, respectively, cadlag functions $[0,T]\to\calS$, which is endowed with the Skorokhod $J_1$-topology.
Write $\C_\calS$ and $\D_\calS$ for the case where $[0,T]$ is replaced by $\R_+$.
Write $X_n\To X$ for convergence in distribution.
A sequence of processes $X_n$ with sample paths in $\D_\calS$ is said to be {\it $C$-tight} if it is tight and every subsequential limit has, with probability 1, sample paths in $\C_\calS$.
Denote by $\io$ the identity map on $\R_+$ defined by $\io(t)=t$ for $t\in\R_+$.

A standard Brownian motion, or SBM for short, is a one-dimensional Brownian motion starting from zero, with zero drift and unit variance.
We abbreviate ``random variable'' and ``independent and identically distributed'' with ``RV'' and, respectively, ``IID''.

\section{The dynamic server population model}\label{sec2}

The queueing model consists of a single queue with multiple servers.
Customers arrive according to a renewal process and are served in the order in which they arrive (i.e., first-come-first-serve).
The service time at each station is exponentially distributed.
Furthermore, when a server becomes idle it waits an exponentially distributed amount of time and then vacations (provided there are still no customers in queue).
Servers spend an exponentially distributed amount of time vacationing before returning to service.

\subsection{Scaling regimes}\label{sec21}

We consider a sequence of queueing networks, indexed by $n\in\N$, that are built on a common probability space $(\Om,\calF,P)$.
For $n\in\N$, let $N^n$ denote the number of servers in the $n$th system, $\la^n>0$ denote the inverse of the mean interarrival times of customers to the system and $\mu^n_{ind}>0$ denote the inverse mean service time.
We use `ind' as a mnemonic for individual service rate.
The parameter $\al\in[\frac{1}{2},1]$ will differentiate the different scaling regimes we consider.
In particular, we assume
\begin{equation*}
	\la^n\sim n,\qquad N^n\sim n^\al,\qquad \mu^n_{ind}\sim n^{1-\al}.
\end{equation*}
Then the overall capacity of the servers, which is given by the product $\mu^n=\mu_{ind}^nN^n$, is of order $n$ and is thus of the same order as the arrival rate $\la^n$.
The three regimes we consider are as follows:
\begin{itemize}
	\item [(i)] $\alpha=1$: HW regime.
	\item [(ii)] $\alpha\in(\frac{1}{2},1)$: near HW regime.
	\item [(iii)] $\alpha=\frac{1}{2}$: NDS regime.
\end{itemize}

\subsection{Customer dynamics}
\label{sec22}

For $n\in\N$ let $Q^n(t)$ denote the number of jobs in the buffer at time $t$.
At any given time, a server can be in three possible states: busy, idle, or vacationing.
The number of servers that are idle and vacationing at time $t$ are denoted by $I^n(t)$ and $V^n(t)$, respectively. 
The number of busy servers is then given by
\begin{equation}\label{27}
	B^n(t):=N^n-I^n(t)-V^n(t),
\end{equation}
and the number of customers in the system at time $t$, denoted by $X^n(t)$, is given by
\begin{equation}
\label{01}
	X^n(t)=Q^n(t)+B^n(t)=Q^n(t)+N^n-I^n(t)-V^n(t).
\end{equation}
The initial conditions $Q^n(0)$, $I^n(0)$ and $V^n(0)$ are $\N_0$-valued RVs represent the number of customers initially in the buffer, the number of servers initially idle and the number of servers initially vacationing, respectively.

Let $\{\IA(l):l\in \mathbb{N}\}$ be strictly positive IID RVs with mean $1$ and variance $C^2_{\IA}>0$, and define
\begin{equation}\label{15}
A^n(t):=\sup\left\{l\geq 0:\sum_{k=1}^l\frac{\IA(k)}{\lambda^n}\leq t\right\}, \qquad t\geq 0\,.
\end{equation}
Then $A^n(t)$ represents the number of customers that arrive in the interval $[0,t]$.
We assume that, as $n\to\iy$,
\begin{align}\label{07}
\hat\lambda^n:=n^{-\frac12}(\la^n-n\la)\to\hat\la,
\end{align}
where $\lambda>0$ and $\hat{\lambda}\in \mathbb{R}$ are constants.
In the $n$th system the server pool consists of $N^n=\lceil n^\al\rceil$ servers.
Each of the servers has IID exponential service times with parameter $\mu_{ind}^n$. 
Recall the overall capacity, given by the product $\mu^n=\mu_{ind}^nN^n$, is of order $n$.
We assume that, as $n\to\iy$,
\begin{equation}\label{08}
\hat\mu^n:=n^{-\frac12}(\mu^n-n\mu)=n^{-\frac12}(\mu_{ind}^nN^n-n\mu)\to\hat\mu,
\end{equation}
where $\mu>0$ and $\hat\mu\in\R$ are constants.
The critical load condition, $\mu=\la$, is assumed throughout this work.
Fix a standard unit Poisson process $S$.
Then the potential service process, denoted by $S^n$, is given by $S^n(t):=S(\mu^nt)$ for $t\ge0$.
The number of departing customers by time $t$, denoted $D^n(t)$, is given by
\begin{equation}
\label{04}
D^n(t)=S\left(\mu^n_{ind}\int_0^tB^n(s)ds\right).
\end{equation}
Denoting by $J^n(t)$ the number of jobs routed to the service pool by time $t$
(not counting the initial number), we also have the following balance equations, namely
\begin{align}
  \label{02}
  Q^n(t)&=Q^n(0)+A^n(t)-J^n(t),\\
  \label{03}
  I^n(t)+V^n(t)&=I^n(0)+V^n(0)+D^n(t)-J^n(t).
\end{align}
A work conservation condition is in force, according to which servers may not be idle when there is work in the queue. 
This can be expressed as
\begin{equation}
  \label{05}
  \text{for every $t\ge0$, $Q^n(t)>0$ implies $I^n(t)=0$.}
\end{equation}
The condition above is in force even for $t=0$,
hence a constraint is implicitly assumed regarding
the initial condition $(Q^n(0),I^n(0),V^n(0))$ alluded to above.

\subsection{Server dynamics}\label{sec23}

Our model for server vacations has positive parameters $\bet^n$ and $\gamma^n$, which are assumed to satisfy, as $n\to\iy$, $\bet^n\to\bet>0$ and $\gam^n\to\gam>0$.
(In \S \ref{sec:multi} we treat a more complicated model in which there are multiple vacationing states.)
The model allows a server to start a vacation only if it is idle. 
It can be described as follows: an exponential clock operating at rate $\bet^n$ is started when the server becomes idle, and it if ticks before the server is busy again, it goes on a vacation for a duration that is exponentially distribution with rate $\gam^n$. 
Thanks to the assumed homogeneity of the servers, this mechanism can be modeled by working with the joint idleness process, and the total number of servers vacationing, rather than accounting for each server individually. To this end, let $S_B$ and $S_E$ be two standard Poisson processes, where $B$ and $E$ are used as mnemonics for {\it beginning} and {\it end} of vacation. 
Then the counting processes associated with vacation beginnings and endings are
\begin{equation}\label{10}
V_B^n(t)=S_B\left(\bet^n\int_0^tI^n(s)ds\right),
\qquad
V_E^n(t)=S_E\left(\gam^n\int_0^tV^n(s)ds\right),
\end{equation}
respectively. Thus the number of servers vacationing is given by
\begin{equation}\label{11}
V^n(t)=V^n(0)+V_B^n(t)-V_E^n(t).
\end{equation}
It is assumed that the five objects
$(Q^n(0),I^n(0),V^n(0))$, $\IA(\cdot)$, $S$, $S_B$ and $S_E$
are mutually independent.

\subsection{Statement of main result}\label{sec24}

We can now state the main result on the model in the exponential vacation case,
characterizing the scaling limit of $(X^n,V^n)$. Define the diffusion scaled process
\begin{equation}\label{16}
	\hat X^n=\frac{X^n-N^n}{\sqrt n},
\end{equation}
where we recall that $N^n=\lceil n^\alpha\rceil$.
Normalize the vacation population size process with scaling specific to $\al$, namely
\begin{equation}\label{17}
	\tilde V^n=\frac{V^n}{n^{\al-\frac12}}.
\end{equation}
Set 
	$$b:=\hat\la-\hat\mu,\qquad\sig^2:=\mu(C_{\IA}^2+1).$$
The statement of the result uses the following terminology. 
Given an $\R_+$-valued process $X=\{X(t),t\ge0\}$, we say that a $\R_+$-valued process $L=\{L(t),t\ge0\}$ is {\it a boundary term for $X$ at zero} if a.s., 
\begin{itemize}
	\item [(i)] $L(0)=0$,
	\item [(ii)] the sample paths of $L$ are non-decreasing, and 
	\item [(iii)] $L$ can only increase when $X$ is zero, i.e., $\int_{[0,\iy)}X(t)dL(t)=0$.
\end{itemize}

\begin{theorem}
  \label{th1}
Fix $\al\in[\frac12,1]$.
Assume that the rescaled initial conditions of $X^n$ and $V^n$ converge, namely that $(\hat X^n(0),\tilde V^n(0))\To(X_0,V_0)$ as $n\to\infty$.
In the case $\al\in[\frac12,1)$, assume also that $X_0\ge0$ a.s.
Then $(\hat X^n,\tilde V^n)\To(X,V)$ as $n\to\infty$, where the pair $(X,V)$ satisfy coupled equations that depend on $\al$ as follows.
\begin{itemize}
\item [(i)] (HW regime) In the case $\al=1$, the pair $(X,V)$ takes values in $\R\times\R_+$ and forms a solution to the SDE-ODE system
\begin{equation}\label{21}
\begin{cases}
\ds
X(t)=X_0+\int_0^t[b+\mu\max(-X(s),V(s))]ds+\sig W(t),
\\ \\
\ds
V(t)=V_0+\int_0^t[\bet(X(s)+V(s))^--\gam V(s)]ds,
\end{cases}
\end{equation}
where $W$ is an SBM, independent of $(X_0,V_0)$.
\item [(ii)] (near-HW regime) In the case $\al\in(\frac12,1)$, the pair $(X,V)$ takes values in $\R_+\times\R_+$ and forms a solution to the SDER-ODE system
\begin{equation}\label{22}
\begin{cases}
\ds
X(t)=X_0+\int_0^t[b+\mu V(s)]ds+\sig W(t)+L(t),
\\ \\
\ds
V(t)=V_0-\gam \int_0^t V(s)ds +\bet \mu^{-1} L(t),
\end{cases}
\end{equation}
where $L$ is a boundary term for $X$ at zero, and $W$ is an SBM, independent of $(X_0,V_0)$.
\item [(iii)] (NDS regime) In the case $\al=\frac12$, the pair $(X,V)$ takes values in $\R_+\times\Z_+$ and forms a solution to the system
\begin{equation}\label{23}
\begin{cases}
\ds
X(t)=X_0+\int_0^t[b+\mu V(s)]ds+\sig W(t)+L(t),
\\ \\
\ds
V(t)=V_0-S_E\left(\gam\int_0^tV(s)ds\right)+S_B(\bet\mu^{-1} L(t)),
\end{cases}
\end{equation}
where $L$ is a boundary term for $X$ at zero, $W$ is an SBM, $S_B$ and $S_E$ are standard Poisson processes, and $W$, $S_B$, $S_E$ and $(X_0,V_0)$ are mutually independent.
\end{itemize}
In case (i) (resp., (ii), (iii)), the system of equations \eqref{21} (resp., \eqref{22}, \eqref{23}) uniquely characterizes the law of the pair $(X,V)$.
\end{theorem}

\begin{remark}\label{rem1}
{\bf (a)}
It is argued in Appendix \ref{sec:a2} that for $\al\in[0,\frac12)$
(the conventional and near-conventional regimes) the unnormalized process
$V^n$ simply vanished in the limit $n\to\iy$. Hence there can be no rescaling under which
the pair of processes remains coupled.
Thus the meaningful regimes for the model studied in this paper are only $\alpha\in[\frac{1}{2},1]$.
\\
{\bf(b)}
When we take $V\equiv 0$ in \eqref{21} we recover the well known SDE studied in \cite{halfin1981heavy}.
Similarly when we take $V\equiv 0$ in \eqref{22} or \eqref{23} we recover the one-dimensional RBM obtained in \cite{atar2012diffusion} when the abandonment rate is set to zero.
\end{remark}

\begin{remark}
Given a solution $(X,V)$ to either \eqref{22} or \eqref{23}, let $L$ be a boundary term for $X$ at zero.
Define the process $\xi=\{\xi(t),t\ge0\}$ by
	$$\xi(t)=X_0+\int_0^t[b+\mu\max(-X(s),V(s))]ds+\sigma W(t).$$
Then the pair $(X,L)$ is a solution to the well known one-dimensional Skorokhod problem for $Z$, and is explicitly given by $(X,L)=\Gam(Z)$, where $\Gam=(\Gam_1,\Gam_2)$ is the one-dimensional Skorokhod map (see Appendix \ref{SP}).
\end{remark}

\section{Proof of main result}\label{sec3}

\subsection{Useful identities and preparatory lemmas}

We first define some related scaled processes.
Namely, let
\begin{align}
	\label{19}
	&\hat A^n(t)=\frac{A^n(t)-\la^nt}{\sqrt n},&&\hat S^n(t)=\frac{S(nt)-nt}{\sqrt n},\\
	\label{24}
	&\hat Q^n(t)=\frac{Q^n(t)}{\sqrt n},&&\tilde I^n(t)=\frac{I^n(t)}{n^{\al-\frac12}}.
\end{align}
The various convergence results in Theorem \ref{th1} are largely based on the fact that centered renewal processes (with finite second moments) satisfy a functional central limit theorem.
In particular, Theorem 14.1 of \cite{Bill} and the mutual independence of the processes $\hat A^n$ and $\hat S^n$ imply that these processes jointly converge to processes $\hat A$ and $\hat S$, which are mutually independent driftless BMs with diffusion coefficients $\sqrt{\la}C_{\IA}$ and, respectively, $1$.

We next develop several identities satisfied by the scaled processes. 
Using first \eqref{01}, \eqref{02} and \eqref{03}, and then \eqref{04}, we obtain
\begin{align*}
	\hat X^n(t)&=\hat X^n(0)+\frac{A^n(t)-D^n(t)}{\sqrt n}\\
	&=\hat X^n(0)+\sigma W^n(t)+n^{-\frac12}\la^nt-n^{-\frac12}\mu^n_{ind}\int_0^tB^n(s)ds,
\end{align*}
where
\begin{equation}
  \label{06}
  W^n(t)=\frac{\hat A^n(t)-\hat S^n\left(n^{-1}\mu^n_{ind}\int_0^t B^n(s)ds\right)}{\sig}.
\end{equation}
Thus by \eqref{27}, \eqref{07}, \eqref{08} and the critical load condition $\la=\mu$, we have, denoting $b^n=\hat\la^n-\hat\mu^n$,
\begin{equation}\label{09}
	\hat X^n(t)=\hat X^n(0)+\sig W^n(t)+b^nt+n^{-1}\mu^n\int_0^t(\tilde I^n(s)+\tilde V^n(s))ds.
\end{equation}
By \eqref{01},
\[
\hat Q^n=\hat X^n+n^{-\frac12}(I^n+V^n)
=\hat X^n+n^{\al-1}\tilde I^n+n^{\al-1}\tilde V^n.
\]
In view of the non-idling condition \eqref{05}, this gives the identities
\begin{equation}\label{14}
(\hat X^n+n^{\al-1}\tilde V^n)^+=\hat Q^n, \qquad (n^{1-\al}\hat X^n+\tilde V^n)^-=\tilde I^n.
\end{equation}

Define the process $Y^n=\{Y^n(t),t\ge0\}$ by
\begin{equation}\label{Yn}
Y^n:=(\hat X^n,\tilde V^n),
\end{equation}
and for a constant $c_0>0$, define the stopping time $\tau^n(c_0)$ by
\begin{equation}\label{taun}
\tau^n(c_0):=\inf\{t:\|Y^n(t)\|\ge c_0\}.
\end{equation}

\begin{lemma}\label{lem:Wn}
Let $\al\in[\frac12,1]$.
Then the processes $W^n$ are $C$-tight.
\end{lemma}

\noi{\bf Proof.}
Recalling that $B^n(t)\le N^n$ by definition \eqref{27}, and using the finiteness of
the constant $c=\sup_nn^{-1}\mu_{ind}^nN^n<\infty$, we have $n^{-1}\mu_{ind}^n\int_0^t B^n(s)ds\le ct$ for all $t\ge0$ and $n\in\N$.
Since the centered renewal processes $\hat A^n$ and $\hat S^n$ are $C$-tight, it follows that the processes $W^n$ are also $C$-tight.
\qed

\skp

The following relations will be useful for cases (i) and (ii), i.e., for $\al\in(\frac12,1]$.
For such $\al$ dividing by $n^{\al-\frac12}$ in \eqref{10} and \eqref{11} yields the relation
\begin{equation}\label{12}
\tilde V^n(t)=\tilde V^n(0)+\bet^n\int_0^t\tilde I^n(s)ds-\gam^n\int_0^t\tilde V^n(s)ds+e^n(t),
\end{equation}
where, with
\begin{equation}
\label{18}
e^n_B(u)=n^{-\al+\frac12}\left[S_B(n^{\al-\frac12}u)-n^{\al-\frac12}u\right],
\qquad
e^n_E(u)=n^{-\al+\frac12}\left[S_E(n^{\al-\frac12}u)-n^{\al-\frac12}u\right],
\end{equation}
we have denoted
\begin{equation}
  \label{13}
  e^n(t)=e_B^n\left(\bet^n\int_0^t\tilde I^n(s)ds\right)-e_E^n\left(\gam^n\int_0^t\tilde V^n(s)ds\right).
\end{equation}

\begin{lemma}\label{lem:en}
Let $\al\in(\frac12,1]$.
Then for all $T<\infty$, $\|e_B^n\|_T+\|e_E^n\|_T\To0$ as $n\to\infty$.
\end{lemma}

\noi{\bf Proof.}
Let $T<\infty$.
The convergence $\|e_B^n\|_T+\|e_E^n\|_T\To0$ follows from definition \eqref{18} and the functional law of large numbers.
\qed

\skp

Next, we record relations that will useful for cases (ii) and (iii), i.e., for $\al\in[\frac12,1)$.
Let 
\begin{equation}\label{eq:tildeX}
	\tilde X^n(t)=(\hat X^n(t))^+\quad\text{and}\quad e^n_X(t)=(\hat X^n(t))^-=\tilde X^n(t)-\hat X^n(t).
\end{equation}
Then $\tilde X^n$ is non-negative and by \eqref{09}, $\tilde X^n=\xi^n+L^n$, where
\begin{equation}\label{34}
	\xi^n(t)=\hat X^n(0)+\int_0^t[b^n+n^{-1}\mu^n\tilde V^n(s)]ds+\sig W^n(t)+e^n_X(t).
\end{equation}
and
\begin{equation}
	\label{28}
	L^n(t)=n^{-1}\mu^n\int_0^t\tilde I^n(s)ds.
\end{equation}
Furthermore, by \eqref{14}, $\tilde I^n(t)>0$ implies $\hat X^n(t)<0$. 
As a result, by \eqref{28} and \eqref{eq:tildeX}, we see that $L^n$ is non-decreasing and $L^n$ can only increase when $\tilde X^n$ is zero, i.e., $\int_0^\infty 1_{\{\tilde X^n(t)>0\}}dL^n(t)=0$.
Consequently, $(\tilde X^n,L^n)$ is the solution to the one-dimensional Skorokhod problem for $\xi^n$ (see Appendix \ref{SP}), and so the pair can be expressed in terms of the one-dimensional Skorokhod map, as follows,
\begin{equation}
\label{30}
(\tilde X^n,L^n)=\Gam(\xi^n).
\end{equation}

\begin{lemma}\label{lem:eXn}
Suppose $\al\in[\frac12,1)$.
Then $e_X^n\To0$.
\end{lemma}

\noi{\bf Proof.}
Fix $\al\in[\frac12,1)$ and $T<\infty$.
By assumption, the weak limit $X_0$ of $\hat X^n(0)$ is non-negative, so it suffices to show that for any $\eps>0$, $P(\Om^{n,\eps})\to0$, where
\[
\Om^{n,\eps}=\left\{\exists\; 0\le s^n< t^n\le T:\hat X^n(s^n)>-2\eps,\hat X^n(t^n)<-3\eps, \sup_{t\in[s^n,t^n]}\hat X^n(t)<-\eps\right\}.
\]
By \eqref{09}, on the event $\Om^{n,\eps}$, there exists $0\le s^n< t^n\le T$ such that
\[
-\eps>\hat X^n(t^n)-\hat X(s^n)\ge-\sig w_T(W^n,\del^n)-c\del^n+n^{-1}\mu^n\int_{s^n}^{t^n}(\tilde I^n(u)+\tilde V^n(u))du,
\]
where $c=\sup_n\|b^n\|<\iy$ and $\del^n=t^n-s^n$.
It follows from \eqref{14} and the fact that $\hat X^n(t)<-\eps$ on the interval $[s^n,t^n]$, that $\tilde I^n+\tilde V^n>n^{1-\al}\eps$ on the interval.
As a result, for all sufficiently large $n$, on $\Om^{n,\eps}$,
	$$n^{-1}\mu^n\int_{s^n}^{t^n}(\tilde I^n(u)+\tilde V^n(u))du\ge c_1\eps n^{1-\al}\del^n,$$ 
where $c_1=\inf_nn^{\al-1}\mu^n_{ind}=\inf_n n^{-1}\mu>0$.
Thus
\[
P(\Om^{n,\eps})\le P(\sig w_T(W^n,\del^n)+c\del^n\ge\eps+c_1\eps n^{1-\al}\del^n).
\]
Fix $\al'\in(0,1-\al)$. 
Separating the cases $\del^n<n^{-\al'}$ and $\del^n\ge n^{-\al'}$, we obtain that, for all sufficiently large $n$,
\[
P(\Om^{n,\eps})\le P(\sig w_T(W^n,n^{-\al'})+cn^{-\al'}\ge\eps)
+P(2\sig\|W^n\|_T+cT\ge c_1\eps n^{1-\al-\al'}).
\]
Both terms on the RHS converge to zero by the $C$-tightness of $W^n$ shown in Lemma \ref{lem:Wn}.
Since $\eps$ and $T$ are arbitrary, this shows that $e_X^n\To0$.
\qed

\subsection{Proof of Theorem \ref{th1}}

We can now prove our main scaling limit result.
Throughout the proof, $C^n$ denotes a generic sequence of constants that satisfy $C^n\to1$, where by the term `generic' we mean that the values the sequence takes may vary from one line to another. 
In addition, the symbol $c$ denotes a generic positive constant (that, in particular, does not depend on $n$).

\skp

\noi{\bf Proof of Theorem \ref{th1}.}
The three regimes are treated separately.
For each regime, our approach is as follows: (a) prove uniqueness in law of solutions to the limiting equations [i.e., \eqref{21}, \eqref{22} or \eqref{23}], (b) express the equations for $(\hat X^n,\tilde V^n)$ in a form that resembles the limiting equations, (c) prove, for each $T$, tightness of $\|Y^n\|_T$, which implies $C$-tightness of the relevant processes, and (d) show that along any convergent subsequent, the limiting processes satisfy the limiting equations (for which uniqueness in law holds).

\skp

{\it (i) The case $\al=1$.}
We first establish uniqueness in law of solutions to the system of equations \eqref{21}.
Observe that \eqref{21} can be viewed as a degenerate SDE with Lipschitz drift and diffusion coefficients (the latter is constant).
For such as SDE, pathwise uniqueness of solutions holds, and consequently so does uniqueness in law.
For the former see Ch.\ V, Theorem 7 of \cite{protter}.

We now write relations for the $(\hat X^n,\tilde V^n)$ that closely resemble the limiting system \eqref{21}.
Setting $\al=1$ in \eqref{14} shows that $\tilde I^n+\tilde V^n=\max(-\hat X^n,\tilde V^n)$, so
\begin{equation}\label{30a}
	\tilde I^n(t)+\tilde V^n(t)\le2\|Y^n(t)\|.
\end{equation}
In addition, substituting the relation into \eqref{09}, and recalling relation \eqref{12} for $\tilde V^n$, yields the system of equations
\begin{equation}\label{26}
\begin{cases}
\ds
\hat X^n(t)=\hat X^n(0)+C^n\mu\int_0^t[b^n+\max(-\hat X^n(s),\tilde V^n(s))]ds+\sig W^n(t),
\\ \\
\ds
\tilde V^n(t)=\tilde V^n(0)+\int_0^t[\bet^n(\hat X^n(s)+\tilde V^n(s))^--\gam^n\tilde V^n(s)]ds+e^n(t).
\end{cases}
\end{equation}

Next, for fixed $T<\infty$, we prove tightness of $\|Y^n\|_T$.
By \eqref{26} and the boundedness of $b^n$, $\beta^n$ and $\gamma^n$, we have
\begin{equation}\label{25a}
	\|Y^n\|_t\le\|Y^n(0)\|+\sig\|W^n\|_t+c_1t+c_1\int_0^t\|Y^n\|_sds+\|e^n\|_t,
\qquad t\ge0.
\end{equation}
for some fixed constant $c_1$ that does not depend on $n$.
Appealing to Gronwall's lemma shows that, for $t\ge0$,
\begin{equation}\label{25b}
	\|Y^n\|_t\le(\|Y^n(0)\|+\sigma\|W^n\|_t+c_1t+\|e^n\|_t)\exp(c_1t).
\end{equation}
Recall that the RVs $\|Y^n(0)\|+\sig\|W^n\|_T$ form a tight sequence by the assumed convergence of the initial conditions and the $C$-tightness of $W^n$ shown in Lemma \ref{lem:Wn}. 
Given $\eps>0$ let $K$ be sufficiently large so that $\limsup_nP(\|Y^n(0)\|+\sig\|W^n\|_T>K)<\eps$.
Choose $c_0>(K+c_1T+1)\exp(c_1T)$ and let $\tau^n=\tau^n(c_0)$ be defined as in \eqref{taun}.
In this case, by \eqref{30a} and the definition of $\tau^n$, $\tilde I^n(t\w\tau^n)+\tilde V^n(t\w\tau^n)\le 2c_0$ for all $t\ge0$, so, in view of definition \eqref{13} for $e^n$ and Lemma \ref{lem:en}, we have $\|e^n\|_{t\w\tau^n}\le\|e^n_B\|_{2c_0\beta^n}+\|e^n_E\|_{2c_0\gamma^n}\To0$ as $n\to\infty$.
Therefore, by our choice of $K$, 
	$$\limsup_{n\to\infty}P(\|Y^n(0)\|+\sig\|W^n\|_T+\|e^n\|_{T\w\tau^n}>K+1)<\eps.$$
From \eqref{25b} and our choice of $c_0$ we see that on the event $\|Y^n(0)\|+\sig\|W^n\|_T+\|e^n\|_{T\w\tau^n}\le K+1$, we have $\|Y^n\|_{T\w\tau^n}<c_0$. 
Hence by definition \eqref{taun}, $\tau^n>T$ on that event.
As a result, $\limsup_nP(\|Y^n\|_T>c_0)<\eps$.
Since $\eps$ is arbitrary, this shows that $\|Y^n\|_T$ forms a tight sequence of RVs.

Having shown tightness of the RVs $\|Y^n\|_T$, we now establish that $(W^n,e^n)\To(W,0)$.
Using \eqref{27}, \eqref{30a} and the fact that $N^n=n$, we have $1-2n^{-\frac12}\|Y^n\|_T\le n^{-1}B^n(t)\le 1$ for all $0\le t\le T$ and $n\in\N$.
As a result, $\|n^{-1}B^n-1\|_T\To0$ as $n\to\iy$. 
Using this, along with the convergence $\mu^n_{ind}\to\mu$ as $n\to\infty$, in the definition \eqref{06} for $W^n$ shows that $W^n(\cdot)\To\sig^{-1}(\hat A(\cdot)-\hat S(\mu\cdot))$.
Note that the latter process is a driftless BM with diffusion coeffient given by $\sig^{-1}(\la C_{\IA}^2+\mu)^{1/2}=\sig^{-1}\la^{1/2}( C_{\IA}^2+1)^{1/2}=1$.
Hence the limit is equal in distribution to the SBM $W$.
Next, observe that the RVs $\int_0^{T}\tilde I^n(s)ds$ and $\int_0^{T}\tilde V^n(s)ds$, which appear as arguments of $e^n_B$ and $e^n_E$ in expression \eqref{13} for $e^n(T)$, form tight sequences of RVs, as follows from \eqref{30a}. 
This, along with Lemma \ref{lem:en}, yields $\|e^n\|_T\To0$.

Now, by \eqref{26} and the established tightness results, we see that the sequence $(\hat X^n,\tilde V^n)$ is $C$-tight.
Taking limits in \eqref{26} along any convergent subsequence, and using that $(\hat X^n(0),\tilde V^n(0))\To(X_0,V_0)$ by assumption, shows that any limit $(X,V,W)$ of $(\hat X^n,\tilde V^n,W^n)$ satisfies \eqref{21}, where we have used that $(C^n,b^n,\beta^n,\gamma^n)\to(1,b,\beta,\gamma)$.
Since uniqueness in law for solutions to \eqref{21} holds, this completes the proof in the case $\alpha=1$.

\skp

{\it (ii) The case $\al\in(\frac12,1)$}.
As in the previous case we first establish uniqueness in law of solutions to the system of equations \eqref{22}, which can be viewed as a degenerate SDER on the domain $\R_+\times\R$ (the non-negativity of the initial condition $V_0$ and the non-decreasing property of $L$ imply that $V(t)\ge0$ for all $t\ge0$; hence there is no necessity to consider the SDER on the smaller domain $\R_+^2$).
The reflection vector field takes the constant value $(1,\beta\mu^{-1})$ on the boundary $\{0\}\times\R$. Theorem 4.3 of \cite{lio-szn} covers such an SDE with reflection on a bounded domain and provides pathwise uniqueness. 
A standard localization argument yields pathwise uniqueness for the unbounded domain at hand.
This shows that uniqueness in law holds for solutions of \eqref{22}.

Next, we write relations for the pair $(\hat X^n,\tilde V^n)$ that closely resemble the limiting system \eqref{22}.
By \eqref{09} and \eqref{12}, for $t\ge0$,
\begin{equation}\label{29}
\begin{cases}
\ds
\hat X^n(t)=\hat X^n(0)+\int_0^t[b^n+C^n\mu\tilde V^n(s)]ds+\sig W^n(t)+L^n(t),
\\ \\
\ds
\tilde V^n(t)=\tilde V^n(0)-\gam^n \int_0^t\tilde V^n(s)ds +C^n\bet^n \mu^{-1} L^n(t)+e^n(t),
\end{cases}
\end{equation}
where $L^n$ is defined as in \eqref{28} and we recall that relation \eqref{30} holds.

We now turn to the proof that for fixed $T<\infty$ the RVs $\|Y^n\|_T$ are tight.
Fix $T<\infty$.
By \eqref{34} for $\xi^n$, we have, for all $t\ge0$,
\begin{equation}\label{29b}
	\|\xi^n\|_t\le\|\hat X^n(0)\|+ct+c\int_0^t\|Y^n\|_sds+c\|W^n\|_t+\|e_X^n\|_t.
\end{equation}
By \eqref{30} and the Lipschitz continuity of $\Gam_2$ (see Proposition \ref{prop:sp}), $L^n(t)\le\|\xi^n\|_t$ for all $t\ge0$.
Thus, using \eqref{29}, we get, for any $t\ge0$, 
\begin{align*}
	\|Y^n\|_t
	&\le c_2\left(Z^n+\|e^n\|_t+\|e^n_X\|_t\right)+c_2\int_0^t\|Y^n\|_sds,
\end{align*}
where $c_2$ is a fixed constant that does not depend on $n$ and
\begin{equation}
	\label{eq:Zn}Z^n:=\|Y^n(0)\|+\|W^n\|_T+T.
\end{equation}
Hence by Gronwall's lemma, for $t\in[0,T]$,
\begin{equation}\label{33}
	\|Y^n\|_t\le c_2\left(Z^n+\|e^n_X\|_t+\|e^n\|_t\right)\exp\left(c_2t\right).
\end{equation}
By the assumed tightness of the initial conditions and $C$-tightness of $W^n$ shown in Lemma \ref{lem:Wn}, the RVs $Z^n$ are tight.
Given $\eps>0$ let $K$ be sufficiently large so that $\limsup_nP(Z^n\ge K)<\eps$. 
Choose $c_0>c_2(K+1)\exp(c_2T)$ and let $\tau^n=\tau^n(c_0)$.
Since $\tilde V^n(T\w\tau^n)\le\|Y^n(T\w\tau^n)\|\le c_0$, it holds that the RVs $\int_0^{T\w\tau^n}\tilde V^n(s)ds$ are tight. 
In addition, by \eqref{28}, the fact that $L^n(t)\le\|\xi^n\|_t$ for all $t\ge0$, and \eqref{29b},
	$$n^{-1}\mu^n\int_0^{T\w\tau^n}\tilde I^n(s)ds=L^n(T\w\tau^n)\le\|\hat X^n(0)\|+c(1+c_0)T+c\|W^n\|_{T\w\tau^n}+\|e_X^n\|_{T\w\tau^n}.$$
Since $\hat X^n(0)$ are tight by assumption, $W^n$ are $C$-tight by Lemma \ref{lem:Wn} and $\|e_X^n\|_{T\w\tau^n}\To0$ by Lemma \ref{lem:eXn}, it follows that the RVs $\int_0^{T\w\tau^n}\tilde I^n(s)ds$ are tight.
In view of definition \eqref{13} for $e^n$ and Lemma \ref{lem:en}, we see that $\|e^n\|_{T\w\tau^n}\To0$ as $n\to\infty$.
Thus, by our choice of $K$,
	$$\limsup_{n\to\infty}P(Z^n+\|e_X^n\|_{T\w\tau^n}+\|e^n\|_{T\w\tau^n}>K+1)<\eps.$$
From \eqref{33} and our choice of $c_0$ we see that on the event $Z^n+\|e_X^n\|_{T\w\tau^n}+\|e^n\|_{T\w\tau^n}\le K+1$, we have $\|Y^n\|_{T\w\tau^n}<c_0$.
Hence by definition \eqref{taun}, $\tau^n>T$ on that event.
As a result, $\limsup_nP(\|Y^n\|_T>c_0)<\eps$.
Since $\eps$ is arbitrary, this shows that $\|Y^n\|_T$ forms a tight sequence of RVs.

Having shown tightness of the RVs $\|Y^n\|_T$, it follows that the RVs $\int_0^T\tilde V^n(s)ds$ are tight.
Hence, by \eqref{29a}, the RVs $\xi^n$ are $C$-tight, and since $(\tilde X^n,L^n)=\Gam(\xi^n)$, it follows from Proposition \ref{prop:sp} that the RVs $(\tilde X^n,L^n)$ are also $C$-tight.
Consequently, since $e_X^n\To0$ (Lemma \ref{lem:eXn}), the RVs $\hat X^n$ are also $C$-tight.
Using that $L^n(t)=n^{-1}\mu^n\int_0^t\tilde I^n(s)ds$ and that $n^{-1}\mu^n\to \mu$ as $n\to\infty$, we see that the RVs $\int_0^t\tilde I^n(s)ds$ are tight.
Thus, the definition \eqref{13} for $e^n$ gives $e^n\To0$ as $n\to\infty$.
Finally, using \eqref{27} and \eqref{28} yields 
	$$N^n-c(\|Y^n\|_t+L^n(t))\le \int_0^t B^n(s)ds\le N^n.$$
As a result, $\|n^{-\al}B^n-1\|_T\To0$ as $n\to\infty$. Using this, along with the convergence $n^{-1+\al}\mu_{ind}^n\to \mu$, in the definition \eqref{06} for $W^n$ shows that $W^n\To W$. 
Taking limits in \eqref{29} and \eqref{30} along any convergent subsequence, and using that $(\hat X(0),\tilde V(0))\To(X_0,V_0)$ by assumption, the Skorokhod representation theorem, the continuity of $\Gam$ and the convergence $e_X^n\To0$, shows that any limit $(X,X,V,W,L,\xi)$ of $(\hat X^n,\tilde X^n,\tilde V^n,W^n,L^n,\xi^n)$ satisfies the equation \eqref{22}, and $(X,L)=\Gam(\xi)$, which, by definition, implies that $L$ is a boundary term for $X$. 
Since uniqueness in law holds for solutions of \eqref{22}, this completes the proof in the case $\alpha\in(\frac12,1)$.

\skp

{\it (iii) The case $\al=\frac12$}.
First note that thanks to the non-negativity of $X_0$, the continuity of the second and third terms in the first equation in \eqref{23} and the oscillation inequality for the SM (Proposition \ref{prop:sp}), the processes $X$ and $L$ have continuous sample paths a.s..
Also, by the second equation in \eqref{23}, $V$ has piecewise constant, right-continuous sample paths.
Uniqueness in law of solutions to the system of the equations \eqref{23} follows from pathwise uniqueness, which we now establish. 
That is, given $W$, $S_B$, $S_E$ and $(X_0,V_0)$, then any two solutions of \eqref{23} are equal a.s.
To this end, let $\Sig=(X,L,V)$ and $\Sig'=(X',L',V')$ be two solutions, let $\tau=\inf\{t:\Sig(t)\ne \Sig'(t)\}$, and consider the event $\tau<\iy$.
By definition, $(X(t),V(t))=(X'(t),V'(t))$ for $t\in[0,\tau)$. 
We first show that $V(\tau)=V'(\tau)$.
To see this must hold, suppose $V$ has a jump of size $-1$ at $\tau$.
Then necessarily $S_E(u)=S_E(u-)+1$, where $u=\gamma\int_0^\tau V(s)ds$.
However, $V=V'$ on $[0,\tau)$, this implies that $V'$ also has a jump of size $-1$ at $\tau$.
A similar argument holds for a jump of size $+1$.
This shows that $V$ and $V'$ agree on $[0,\tau]$ whenever $\tau<\iy$.
Since $V$ is piecewise constant with right-continuous sample paths and $V(\tau)=V'(\tau)$, there exists $\eps>0$ (depending on the sample path) such that $V(t)=V'(t)$ for all $t\in[0,\tau+\eps)$.
It follows from the expression for $(X,L)$ and $(X',L')$ in terms of the one-dimensional Skorokhod map that they are also equal on $[0,\tau+\eps)$, thus contradicting the definition of $\tau$.
With this contradiction thus obtained, we must have $\tau=\infty$ a.s., so pathwise uniqueness holds.

Next, we write relations for the $(\hat X^n,\tilde V^n)$ that closely resemble the limiting system \eqref{23}.
By \eqref{09} and \eqref{10}--\eqref{11},
\begin{equation}\label{29a}
\begin{cases}
\ds
\hat X^n(t)=\hat X^n(0)+\int_0^t[b^n+C^n\mu V^n(s)]ds+\sig W^n(t)+L^n(t),
\\ \\
\ds
V^n(t)=V^n(0)-S_E\left(\gam^n \int_0^t V^n(s)ds\right)+S_B\left(C^n\bet^n \mu^{-1} L^n(t)\right),
\end{cases}
\end{equation}
where $L^n$ is defined as in \eqref{28}.
In addition we recall that \eqref{30} holds.

We now turn to the proof of tightness of the RVs $\|Y^n\|_T$.
The main difference between this case and the case $\al\in(\frac12,1)$ is the treatment of the equation that governs $V^n$, where we note that $\tilde V^n=V^n$ in this case.
We argue as follows.
By \eqref{34}, for all $t\ge0$,
\begin{equation}\label{eq:xinbound}
	\|\xi^n\|_t\le c_4M^n+c_4\int_0^t V^n(s)ds,
\end{equation}
where $M^n:=\|\hat X^n(0)\|+T+\|W^n\|_T+\|e^n_X\|_T$ and $c_4\ge1$ is a suitable constant that does not depend on $n$. 
Hence by \eqref{29a} and the fact that $L^n(t)\le\|\xi^n\|_t$,
\begin{align}\label{36}
	\|V^n\|_t&\le V^n(0)+S_B(c\|\xi^n\|_t)\le V^n(0)+S_B\left(c_4 M^n+c_4\int_0^tV^n(s)ds\right),
\end{align}
where we have chosen $c_4$ to be possibly larger.
By the convergence of the initial conditions $(V^n(0),\hat X^n(0))\Rightarrow (V_0,X_0)$, the $C$-tightness of $W^n$ (Lemma \ref{lem:Wn}) and the convergence $e_X^n\To0$ (Lemma \ref{lem:eXn}), it follows that the RVs $V^n(0)+c_4M^n$ are tight.
Let $\eps>0$ and choose $K<\infty$ sufficiently large such that $P(V^n(0)+c_4M^n\ge K)<\eps$ for all $n$.
In addition, since $S_B$ is a Poisson process, by the functional law of large numbers, $\|k^{-1}S_B(k\cdot)-\iota(\cdot)\|_u\To0$ as $k\to\infty$, for $u=1+2c_4Te^{c_4T}$.
Thus, by choosing $K<\infty$ possibly larger, we can ensure that $P(\Om^n)>1-2\eps$, where
	$$\Om^n:=\left\{V^n(0)+c_4M^n\le K\text{ and }S_B(t)<t+K\text{ for all }t\le K+2c_4KTe^{c_4T}\right\}.$$
Let $c_0:=2K+4c_4KTc^{e_4T}+2Ke^{c_4T}$ and $\tau^n=\tau^n(c_0)$.
Then $c_4M^n+c_4\int_0^tV^n(s)\le K+2c_4KTe^{c_4T}$ for all $t\le T\w\tau^n$.
Thus, on the event $\Om^n$, we have, for all $t\le T$,
	$$\|V^n\|_{t\w\tau^n}\le V^n(0)+c_4M+c_4\int_0^{t\w\tau^n}V^n(s)ds+K\le 2K+c_4\int_0^{t\w\tau^n} V^n(s)ds.$$
Hence, by Gronwall's lemma, on this event, $\|V^n\|_{T\w\tau^n}\le 2Ke^{c_4T}$.
Moreover, by \eqref{30}, the Lipschitz continuity of the SM (Proposition \ref{prop:sp}) and \eqref{eq:xinbound}, on this event we have
	$$\|\hat X^n\|_{T\w\tau^n}\le 2\|\xi^n\|_{T\w\tau^n}\le 2c_4M^n+2c_4\int_0^tV^n(s)\le 2K+4c_4KTe^{c_4T}.$$
Combining the bounds on $V^n$ and $\hat X^n$ gives $\|Y^n\|_{T\w\tau^n}\le 2K+4c_4KTe^{c_4T}+2Ke^{c_4T}$. 
Thus, $\tau^n>T$ on $\Om^n$.
Since $\eps$ is arbitrary, this shows the tightness of the RVs $\|Y^n\|_T$. 

Having shown tightness of the RVs $\|Y^n\|_T$, we can argue exactly as in case (ii) to conclude that $(\hat X^n,\tilde X^n, V^n,L^n,\xi^n)$ are $C$-tight.
Taking limits in \eqref{29a} and \eqref{30} along any convergent subsequence, and using that $(\hat X(0),V(0))\To(X_0,V_0)$ by assumption, the Skorokhod representation theorem, the continuity of $\Gam$ and the convergence $e_X^n\To0$, shows that any limit $(X,X,V,W,L,\xi)$ of $(\hat X^n,\tilde X^n,V^n,W^n,L^n,\xi^n)$ satisfies the equations \eqref{23}, and $(X,L)=\Gam(\xi)$, which, by definition, implies that $L$ is a boundary term for $X$. 
Since uniqueness in law holds for solutions of \eqref{23}, this completes the proof in the case $\alpha=\frac12$.
\qed

\section{Multi-stage vacation model}
\label{sec:multi}

In this section we consider a generalization of the model that allows for multiple vacationing states.
Let $m\ge2$ denote the number of vacationing states.
In the multi-stage model, servers may take the following states: busy, idle, and vacationing state $i$ for $i=1,\dots,m$.
Let $\M=\{1,\dots,m\}$.
At time $t\ge0$ let $I^n(t)$ denote the number of idle servers and $U_i^n(t)$ denote the number of vacationing servers in state $i$, for $i\in\M$.
Let $U^n(t)=(U_1^n(t),\dots,U_m^n(t))$.
Then $V^n(t)=U^n(t)\cdot 1=U_1^n(t)+\cdots+U_m^n(t)$ denotes the total number of vacationing servers and $B^n(t)$, defined as in \eqref{27}, denotes the number of busy servers.

The customer dynamics described in \S \ref{sec22} still hold; in particular, equations \eqref{27}--\eqref{05} hold.
The server dynamics, however, are no longer described by \eqref{10}--\eqref{11}.
For the multi-stage setting, fix vectors $\beta^n,\gamma^n\in\R_+^m$.
Here $\beta_i^n$ denotes the rate at which idling servers transition to vacationing state $i$ and $\gamma_i^n$ denotes the rate at which vacationing servers in state $i$ return to idling, for $i\in\M$.
Fix an $m\times m$ transition rate matrix $R^n=(r_{ij}^n)$ so that $r_{ij}^n$ denotes the rate at which vacationing servers transition from state $i$ to $j$, for $i\ne j\in\M$, and $r_{ii}^n=-\sum_{j\ne i}r_{ij}^n$.
We assume there exist vector $\beta,\gamma\in\R_+^m$ and an $m\times m$ matrix $R$ such that $(\beta^n,\gamma^n,R^n)\to(\beta,\gamma,R)$ as $n\to\infty$.
Let $\M_0=\M\cup\{0\}$ and $S_{ij}$, for $i\ne j\in\M_0$, be independent unit Poisson processes.
For $i\in\mathbb{M}$, define 
\begin{align}\label{ViB}
	U_{i,B}^n(t)&=S_{0i}\left(\beta_i^n\int_0^tI^n(s)ds\right)+\sum_{j\ne i}S_{ji}\left(r_{ji}^n\int_0^tU_j^n(s)ds\right),\\ \label{ViE}
	U_{i,E}^n(t)&=S_{i0}\left(\gamma_i^n\int_0^tU_i^n(s)ds\right)+\sum_{j\ne i}S_{ij}\left(r_{ij}^n\int_0^tU_i^n(s)ds\right).
\end{align}
Then
\begin{equation}\label{Vi}
	U_i^n(t)=U_i^n(0)+U_{i,B}^n(t)-U_{i,E}^n(t).
\end{equation}
It is assumed that the objects $(Q^n(0),I^n(0),U^n(0))$, $IA(\cdot)$, $S$ and $S_{ij}$, $i\neq j\in\M_0$, are mutually independent.
Define $\hat X^n$ as in \eqref{16} and define $\tilde U^n$ by
\begin{equation}\label{17m}
	\tilde U^n=\frac{U^n}{n^{\al-\frac12}}.
\end{equation}
Our main result for the multi-stage vacation model characterizes the limit of the scaled pair $(\hat X^n,\tilde U^n)$.

\subsection{Statement of main result}

Let $G^n=\diag(\gamma^n)$ denote the $m\times m$ diagonal matrix satisfying $G_{ii}^n=\gamma_i^n$ for $i\in\M$.
Recall the definition of a boundary term given prior to Theorem \ref{th1}.
We can now state our main result on limits of the scaled queue-server processes in the case of multi-stage vacations.

\begin{theorem}
	\label{th2}
	Fix $\al\in[\frac12,1]$.
	Assume that the rescaled initial conditions of $X^n$ and $V^n$ converge, namely that $(\hat X^n(0),\tilde U^n(0))\To(X_0,U_0)$.
	In the case $\al\in[\frac12,1)$, assume also that $X_0\ge0$ a.s.
	Then $(\hat X^n,\tilde U^n)\To(X,U)$, where the law of the limit process $(X,U)$ depends on $\al$ and is specified in what follows.
	\begin{itemize}
	\item[(i)] (HW regime) In the case $\al=1$, the pair $(X,U)$ takes values in $\R\times\R_+^m$ and forms a solution to the SDE-ODE system
	\begin{equation}\label{21m}
	\begin{cases}
	\ds
	X(t)=X_0+\int_0^t[b+\mu\max(-X(s),V(s))]ds+\sig W(t),
	\\ \\
	\ds
	U(t)=U_0+\int_0^t[\beta(X(s)+V(s))^- +(R-G)U(s)]ds,	
	\end{cases}
	\end{equation}
	where $V=U\cdot 1$ and $W$ is a SBM, independent of $(X_0,U_0)$.
	
	\item[(ii)] (near-HW regime) In the case $\al\in(\frac12,1)$, the pair $(X,U)$ takes values in $\R_+\times\R_+$ and
	forms a solution to the SDER-ODE system
	\begin{equation}\label{22m}
	\begin{cases}
	\ds
	X(t)=X_0+\int_0^t[b+\mu V(s)]ds+\sig W(t)+L(t),
	\\ \\
	\ds
	U(t)=U_0+\int_0^t (R-G)U(s)ds+\beta\mu^{-1}L(t),
	\end{cases}
	\end{equation}
	where $V=U\cdot 1$, $L$ is a boundary term for $X$ at zero, and $W$ is a SBM, independent of $(X_0,U_0)$.
	
	\item[(iii)] (NDS regime) In the case $\al=\frac12$, the pair $(X,U)$ takes values in $\R_+\times\Z_+$ and
	forms a solution to the system
	\begin{equation}\label{23m}
	\begin{cases}
	\ds
	X(t)=X_0+\int_0^t[b+\mu V(s)]ds+\sig W(t)+L(t),
	\\ \\
	\ds
	U_i(t)=U_{0i}-\sum_{j\in\mathbb{M}\setminus\{i\}}S_{ij}\left(r_{ij}\int_0^tU_i(s)ds\right)+\sum_{j\in\mathbb{M}\setminus\{i\}}S_{ji}\left(r_{ji}\int_0^tU_j(s)ds\right)\\ 
	\ds
	\qquad\qquad+S_{0i}(\beta_i\mu^{-1} L(t))-S_{i0}\left(\gamma_i\int_0^tU_i(s)ds\right),
	\end{cases}
	\end{equation}
	where $V=U\cdot1$, $L$ is a boundary term for $X$ at zero,
	$W$ is a SBM, $S_B$ and $S_E$ are standard Poisson processes,
	and $W$, $S_{ij}$, $i\ne j\in\mathbb{M}_0$, and $(X_0,U_0)$ are mutually independent.
	
	\item[(iv)] In case (i) (resp., (ii), (iii)), the system of equations \eqref{21m} (resp., \eqref{22m}, \eqref{23m}) uniquely characterizes the law of the pair $(X,U)$.
\end{itemize}
\end{theorem}

\begin{remark}
	The equations for $X$ in the multi-stage vacation model are exactly the same as the equations for $X$ in the singe-stage vacation model owing to the fact that the customer dynamics only require information about the total number of servers on vacation.
\end{remark}

\subsection{Proof of Theorem \ref{th2}}

Define $\hat A^n$, $\hat S^n$, $\hat Q^n$ and $\tilde I^n$ as in \eqref{19}--\eqref{24}. 
Then, as stated there, the centered and scaled processes $\hat A^n$ and $\hat S^n$ converge to driftless BMs with diffusion coefficients $\sqrt{\lambda}C_{IA}$ and, respectively, 1.
Define
	$$\tilde V^n=\frac{V^n}{n^{1-\alpha}}=\tilde U^n\cdot 1.$$
Then \eqref{09} holds with $W^n$ defined as in \eqref{06}, and relations \eqref{14} hold.
Define $Y^n$ and $\tau^n(c_0)$, for $c_0>0$, as in \eqref{Yn}--\eqref{taun}.
Then Lemma \ref{lem:Wn} holds by the exact same argument.

We now derive equations for $\tilde U^n$ for cases (i) and (ii), i.e., for $\al\in(\frac12,1]$.
Dividing by $n^{\al-\frac12}$ in \eqref{ViB}--\eqref{Vi}, we obtain
\begin{align}
\label{12m}
\tilde U_i^n(t)&=\tilde U_i^n(0)+\bet_i^n\int_0^t\tilde I^n(s)ds+\sum_{j\in\mathbb{M}\setminus\{i\}}r_{ji}^n\int_0^t\tilde U_j^n(s)ds\\
\notag
&\qquad-\gam_i^n\int_0^t\tilde U_i^n(s)ds-\sum_{j\in\mathbb{M}\setminus\{i\}}r_{ij}^n\int_0^t\tilde U_i^n(s)ds+e_i^n(t),
\end{align}
where, with
\begin{align}
\label{18ma}
e^n_{ij}(u)&=n^{-\al+\frac12}\left[S_{ij}\left(n^{\al-\frac12}u\right)-n^{\al-\frac12}u\right],
\end{align}
for $i\ne j\in\mathbb{M}_0$, we have denoted
\begin{align}
\label{13m}
e_i^n(t)&=e_{0i}^n\left(\bet_i^n\int_0^t\tilde I^n(s)ds\right)+\sum_{j\in\mathbb{M}\setminus\{i\}}e_{ji}^n\left(r_{ji}^n\int_0^t\tilde U_j^n(s)ds\right)\\
\notag
&\qquad-e_{i0}^n\left(\gam_i^n\int_0^t\tilde U_i^n(s)ds\right)-\sum_{j\in\mathbb{M}\setminus\{i\}}e_{ij}^n\left(r_{ij}^n\int_0^t\tilde U_i^n(s)ds\right).
\end{align}
Let
	$$e^n(t)=e_1^n(t)+\cdots+e_m^n(t),\qquad t\ge0.$$

\begin{lemma}\label{lem:enm}
Suppose $\al\in(\frac12,1]$.
Then $e_{ij}^n\To0$ for each $i\ne j\in\mathbb{M}_0$.
\end{lemma}

\noi{\bf Proof.}
This follows from definition \eqref{18ma} and the functional law of large numbers.
\qed

\skp

Next, in cases (ii) and (iii), i.e., for $\al\in[\frac12,1)$, define $\tilde X^n$ and $e_X^n$ as in \eqref{eq:tildeX}, and define $\xi^n$ and $L^n$ as in \eqref{34} and, respectively, \eqref{28}.
Then, by the same argument presented there, \eqref{30} holds, as does Lemma \ref{lem:eXn}.

We can now prove Theorem \ref{th2}.
The convention from the previous section regarding the notation $C^n$ and $c$ is kept.

\skp

\noi{\bf Proof of Theorem \ref{th2}.}
The proof follows a similar structure to the proof of Theorem \ref{th1} and many of the arguments are identical or quite similar.
Here we describe the new aspects of the proof and refer the reader to the proof of Theorem \ref{th1} when the arguments are identical.
As in the proof of Theorem \ref{th1}, the different regimes are treated separately.

\skp

{\it The case $\al=1$.}
Uniqueness in law of solutions to the system of equations \eqref{21m} follows from the fact that the system can be viewed as a degenerate SDE with Lipschitz drift and diffusion coefficients, for which pathwise uniqueness of solutions holds.

Next, we write relations for $(\hat X^n,\tilde U^n)$ that are similar to \eqref{21m}.
As in the single stage setting, equation for $\hat X^n$ in \eqref{26} holds.
After substituting the relation \eqref{14} (with $\al=1$) into equation \eqref{12m} for $\tilde U^n$ with relation and recalling that $R=(r_{ij})$ and $G=\diag(\gamma^n)$, we arrive at the system of equations
\begin{equation}\label{26m}
\begin{cases}
\ds
\hat X^n(t)=\hat X^n(0)+C^n\mu\int_0^t[b^n+\max(-\hat X^n(s),\tilde V^n(s))]ds+\sig W^n(t),
\\ \\
\ds
\tilde U^n(t)=\tilde U^n(0)+\int_0^t[\bet^n(\hat X^n(s)+\tilde V^n(s))^-+(R-G)U^n(s)]ds+e^n(t).
\end{cases}
\end{equation}

We now prove tightness of $\|Y^n\|_T$ for all $T<\infty$.
Let $T<\infty$.
From \eqref{26m} we see that \eqref{25a} holds with a possibly different choice of $c_1$.
The remaining argument that the RVs $\|Y^n\|_T$ are tight follows exactly as in the single stage setting, except that 
	$$\|e^n\|_{t\w\tau^n}\le\sum_{i\ne j\in\M_0}\|e_{ij}^n\|_{2c_0\tilde c}\To0,$$
with $\tilde c=\sup_n\max_{i,j}(\beta_i^n,\gamma_i^n,r_{ij}^n)$, follows from Lemma \ref{lem:enm} instead of Lemma \ref{lem:en}.
Having established tightness of the RVs $\|Y^n\|_T$ for all $T<\infty$, we argue exactly as in the single stage setting that $(W^n,e^n)\To(W,0)$ (again using Lemma \ref{lem:enm} instead of Lemma \ref{lem:en}) and the RVs $\int_0^T\tilde I^n(s)ds$ and $\int_0^T\tilde V^n(s)ds$ form tight sequences of RVs.
Thus, in view of Lemma \ref{lem:enm}, $e^n\To0$.
Taking limits in \eqref{26m} along any convergent subsequence and using that $(\hat X^n(0),\tilde U^n(0))\To(X_0,U_0)$ by assumption, shows that any limit $(X,U,W)$ of $(\hat X^n,\tilde U^n,W^n)$ satisfies \eqref{21m}, where we have used that $(C^n,b^n,\beta^n,\gamma^n,R^n)\to(1,b,\beta,\gamma,R)$.
Since uniqueness in law of solutions to \eqref{21m} holds, this completes the proof.

\skp

{\it The case $\al\in(\frac12,1)$}.
Uniqueness in law of solutions to the system \eqref{22m}, which can be viewed as a degenerate SDER on the domain $\R_+\times\R^m$, with constant reflection vector $(1,\beta_1\mu^{-1},\dots,\beta_m\mu^{-1})$ on the boundary $\{0\}\times\R^m$.
Since pathwise uniqueness of such an SDER holds (again via a localization argument), this implies that uniqueness in law holds for solution of \eqref{22m}.

By \eqref{09} and \eqref{12m},
\begin{equation}\label{29m}
	\begin{cases}
	\ds
	\hat X^n(t)=\hat X^n(0)+\int_0^t[b^n+C^n\mu\tilde V^n(s)]ds+\sig W^n(t)+L^n(t),
	\\ \\
	\ds
	\tilde U^n(t)=\tilde U^n(0)+\int_0^t(R-G)\tilde U^n(s)ds+C^n\beta^n\mu^{-1}L^n(t)+e^n(t),
	\end{cases}
\end{equation}
where $L^n$ is defined as in \eqref{28} and we recall that relation \eqref{30} holds.

The proof that the RVs $\|Y^n\|_T$ for fixed $T<\infty$ follows the same argument as in the single stage setting and, as in that case, implies $C$-tightness of $\xi^n$, $\tilde X^n$, $L^n$ and $\hat X^n$, and the convergence $W^n\To W$.
Taking limits in \eqref{29m} along any convergent subsequence, and using the convergence of the initial conditions, the Skorokhod representation theorem and the continuity of $\Gam$ shows that any limit $(X,X,U,W,L,\xi)$ of $(\hat X^n,\tilde X^n,\tilde U^n,W^n,L^n,\xi^n)$ satisfies equation \eqref{22m} and $(X,L)=\Gam(\xi)$, so $L$ is a boundary term for $X$.
Since uniqueness in law holds for solutions to \eqref{22m}, this completes the proof in the case $\al\in(\frac12,1)$.

\skp

{\it The case $\al=\frac12$}.
Uniqueness in law of solutions to \eqref{23m} follows from pathwise uniqueness, which is shown using an argument analogous to the one in the proof of Theorem \ref{th1}.
By \eqref{09} and \eqref{ViB}--\eqref{Vi},
\begin{equation}\label{29am}
	\begin{cases}
	\ds
	\hat X^n(t)=\hat X^n(0)+\int_0^t[b^n+C^n\mu V^n(s)]ds+\sig W^n(t)+L^n(t),
	\\ \\
	\ds
	U_i^n(t)=U_i^n(0)-\sum_{j\in\mathbb{M}\setminus\{i\}}S_{ij}\left(r_{ij}^n\int_0^tU_i^n(s)ds\right)+\sum_{j\in\mathbb{M}\setminus\{i\}}S_{ji}\left(r_{ji}^n\int_0^tU_j^n(s)ds\right)\\ 
	\ds
	\qquad\qquad+S_{0i}(C^n\beta_i\mu^{-1} L^n(t))-S_{i0}\left(\gamma_i^n\int_0^tU_i^n(s)ds\right),
	\end{cases}
\end{equation}
where $L^n$ is defined as in \eqref{28}.
In addition, we recall that \eqref{30} holds.

Next, we prove tightness of the sequence $\|Y^n\|_T$ for fixed $T<\infty$.
As in the single stage server setting, \eqref{eq:xinbound} holds.
Hence, by \eqref{29am}, for $i\in\M$,
\begin{align}
	\|V^n\|_t&=\|U_1^n+\cdots+U_m^n\|_t\le V^n(0)+\sum_{i\in\M}S_{0i}\left(c_4M^n+c_4\int_0^tU_i^n(s)ds\right),
\end{align}
where we have chosen $c_4$ possibly larger.
As in the single stage setting, the RVs $V^n(0)+c_4M^n$ are tight.
Let $\eps>0$ and choose $K<\infty$ sufficiently large such that $P(V^n(0)+c_4M^n\ge K)<\eps$ for all $n$.
In addition, we can choose $K<\infty$ possibly larger so that $P(\Om^n)>1-2\eps$, where
	$$\Om^n=\left\{V^n(0)+c_4M^n\le K\text{ and }S_{0i}(t)<t+K\text{ for all }t\le K+2c_4KTe^{c_4T}\text{ and }i\in\M\right\}.$$
The remainder of the proof proceeds exactly as in the single-stage vacation setting.
\qed

\section{Heuristic formulas for steady-state performance}\label{sec:heu}

The goal of this section is to quantify the effect of the server vacations on
performance measures that are of particular interest in the HW and the NDS regimes. 
In the HW regime (i.e., when $\al=1$), it was shown in \cite{halfin1981heavy} that the probability of wait (POW) converges to a nondegenerate limit strictly between $0$ and $1$, where as in all other regimes (i.e., when $0\le\al<1$), POW converges to $1$.
Furthermore, an explicit formula for the limit of POW was derived in \cite{halfin1981heavy}.
As argued in \cite{whi-04}, POW is an especially useful performance measure because it requires no scaling by a function on $n$ in this regime.
It is desirable to understand how this limit probability differs in the case of server
vacations of the form that we study.
Potentially, this information could be obtained from the steady state distribution
of the diffusion limit of Theorem \ref{th1}(i). However, one does not expect
an explicit formula for this two-dimensional dynamics.
Instead, we develop in \S \ref{sec41} a formula based on Theorem \ref{th1}(i)
and on a heuristic argument in which we substitute the steady-state mean number of available servers into the explicit formula derived in \cite{halfin1981heavy}.
We provide numerical tests of the accuracy of this heuristic formula.

In the NDS regime (i.e., when $\al=1/2$), the slowdown is a performance measure that is particularly important. 
The slowdown converges (without any rescaling) to a number strictly between $1$ and $+\iy$ in the NDS regime, whereas it converges to either 1 or $+\iy$ in all other regimes (i.e., when $\al\ne1/2$).
A formula for the slowdown in absence of server vacations was given
in \cite{atar2012diffusion}. Again, although an expression for the slowdown, based
on the steady state distribution of the coupled pair of Theorem \ref{th1}(iii), could be developed,
it would not be explicit. Thus instead we develop, in \S\ref{sec42},
a variant of the formula from \cite{atar2012diffusion}, based on Theorem \ref{th1}(iii)
and a heuristic argument similar to the one used for the HW regime, accompanied by numerical tests.

\subsection{Steady-state probability of waiting for service in the HW regime}\label{sec41}

In the original setting of \cite{halfin1981heavy} for the model with no vacations,
the diffusion limit is given by the SDE
\begin{equation}\label{114}
	X(t)=X(0)+\int_0^t[b+\mu(X(s))^-]ds+\sig W(t),
\end{equation}
which is precisely equation \eqref{21} with $V=0$. 
In the case $b<0$, this Markov process has a unique invariant distribution, and so the steady state probability $POW_0:=P(X(t)>0)$
is well defined
as a quantity that does not depend on $t$. We sometimes write this in an informal
fashion as $P(X(\iy)>0)$. It is also proved there that this probability
is the $n\to\iy$ limit of the $n$th system probability of wait.
Throughout, $b<0$ is assumed.
Denote by $\Ph$ the standard normal cumulative distribution function.
Then for the model without vacations, the formula reads
\begin{equation}\label{120}
POW_0=\frac{1}{1+\sqrt{2\pi}b_1\Ph(b_1)\exp(b_1^2/2)}
\qquad
\text{where}
\qquad
b_1=\sqrt{\frac{2}{\mu}}\frac{|b|}{\sig}.
\end{equation}
For the original expression from \cite{halfin1981heavy},
see Th.\ 4 there; a corrected version of this formula
appeared in \cite[eq.\ (1.1)]{whi-04}
(although the original formula is correct in the special case $\sigma=\sqrt{2}$,
corresponding to the M/M/n limit law).
Formula \eqref{120} above follows the corrected version,
albeit with somewhat different notation.

We now go back to the model that accommodates vacations,
in which case the limit is given by \eqref{21}.
For each $n\in\N$, define $Y^n(t)=Q^n(t)-I^n(t)$ for all $t\ge0$, so that, by the non-idling policy,
$Y^n(t)>0$ if and only if $Q^n(t)>0$. On this event, an arriving customer necessarily waits in the queue.
Define the scaled process $\hat Y^n$ by
\[
\hat Y^n(t)=\frac{Y^n(t)}{\sqrt{n}},\qquad t\ge0.
\]
Then (recalling $\al=1$)
\[
\hat Y^n(t)=\hat X^n(t)+\tilde V^n(t),\qquad t\ge0,
\]
and so $\hat Y^n\Rightarrow Y$ as $n\to\infty$,
where $Y(t)=X(t)+V(t)$ for all $t\ge0$.
We are therefore led to study the probability $P(Y(t)>0)$ at steady state (where it is independent of $t$).
We have from \eqref{21}
\begin{equation}\label{21a}
\begin{cases}
\ds
Y(t)=Y_0+\int_0^t[b+(\mu+\beta)(Y(s))^-+(\mu-\gamma)V(s)]ds+\sig W(t),
\\ \\
\ds
V(t)=V_0+\int_0^t[\bet(Y(s))^--\gam V(s)]ds.
\end{cases}
\end{equation}
In general an exact expression is beyond the scope of this work.
Instead, we use heuristics to estimate the probability, as follows.
Denote by $y,v\ge0$ the a.s.\ averages
\[
y=\lim_{t\to\infty}\frac{1}{t}\int_0^t(Y(s))^-ds,\qquad v=\lim_{t\to\infty}\frac{1}{t}\int_0^tV(s)ds.
\]
From \eqref{21a} we obtain that $b+(\mu+\beta)y+(\mu-\gamma)v=0$ and $\beta y=\gamma v$.
Solving for $y$ and $v$ we have
\[
y=\frac{\gamma |b|}{\mu(\gamma+\beta)},\qquad v=\frac{\beta|b|}{\mu(\gamma+\beta)}.
\]
The quantity $v$ is the steady state value of the process $V$. The main heuristic step is now to
substitute $v$ in for $V(s)$ in \eqref{21a}. This yields the (uncoupled) SDE
\begin{equation}\label{21b}
\tilde Y(t)=Y_0+\int_0^t[\tilde b+\tilde\mu(\tilde Y(s))^-]ds+\sigma W(t),
\end{equation}
where we denote
\[
\tilde b=\frac{\gamma(\mu+\beta)}{\mu(\gamma+\beta)}b,
\qquad
\tilde\mu=\mu+\beta.
\]
Then $P(\tilde Y(\infty)>0)\approx P(Y(\infty)>0)$
with equality holding when $\mu=\gamma$ since the SDE for $Y$ in \eqref{36}
becomes uncoupled from the ODE for $V$. In what follows we will use the notation
$POW=P(Y(\iy)>0)$ and $\widetilde{POW}=P(\tilde Y(\iy)>0)$
for the exact and, respectively, approximate performance measure.

For an explicit expression for $\widetilde{POW}$
we only need to notice the similarity of \eqref{21b} to \eqref{114}.
That is, $\tilde Y$ satisfies the same equation that $X$ satisfies in the model with
no vacations, but with different parameters, hence
$P(\tilde Y(\iy)>0)$ can be obtained from \eqref{120} upon modifying the
parameters. This gives
\begin{equation}\label{1a}
\widetilde{POW}=\frac{1}{1+\sqrt{2\pi}b_2\Ph(b_2)\exp(b_2^2/2)}
\qquad
\text{where}
\qquad
b_2=\sqrt{\frac{2}{\tilde\mu}}\frac{|\tilde b|}{\sig}
=\frac{\gamma\sqrt{2(\mu+\beta)}}{\mu(\gamma+\beta)}\frac{|b|}{\sig}.
\end{equation}

The main heuristic step is replacing the stochastic process $V(t)$ by its average.
Hence it is expected that the approximation is good when the time scale
of vacation lengths is long compared to the time scale at which the queue length
fluctuates. We return to this point when discussing the simulation results.

\begin{comment}
The choice of parameters we take for this regime is $\mu=0.1$, $\beta=0.1$, $M=10$, $\sigma=1$, $T=10^5$, $N=10^8$. First, we make sure that the simulation error is small. To ensure this, we simulate the heuristic formula and show that its relative error compared to the already known formula for probability of delay is very small. Then with these simulation parameters we compare probability of delay for $Y(\cdot)$ (calculated using simulation) with the already known formula. The results are shown in  Fig \ref{hw_gamma}. It is clear from Fig \ref{hw_gamma} that the heuristic formula approximates the probability of delay for $Y(\cdot)$ well when $\gamma$ is small and it is bad approximation for large values of $\gamma$. 
\end{comment}

\subsection{Steady-state slowdown in the NDS regime}\label{sec42}

It is well known that the steady state distribution of
\begin{equation}\label{1c}
X(t)=x+bt+\sig W(t)+L(t),
\end{equation}
where $L$ is the boundary term for $X$ at zero and $b<0$,
is exponential with mean $E X(\iy)=\sig^2/(2|b|)$. Based on this and the fact that
the leading term in the expected service time is given by $\mu\sqrt{n}$, the limiting slowdown (SD)
in the NDS regime, for the model with no vacationing servers,
was computed in \cite{atar2012diffusion} to give
\begin{equation}\label{1d}
SD_0=1+EX(\iy)=1+\frac{\sig^2}{2|b|}.
\end{equation}
For the model with vacations, the relevant quantity for similar considerations is
\[
SD=1+EX(\iy),
\]
where $(X,V)$ is a solution to \eqref{23}, and, throughout, $b<0$ is assumed.
As before, we do not aim at an exact calculation because an explicit expression
is not expected; however, again one can proceed via a heuristic argument.
Define
\[
\ell=\lim_{t\to\infty}\frac{L(t)}{t},\qquad v=\lim_{t\to\infty}\frac{1}{t}\int_0^tV(s)ds.
\]
Then by \eqref{23}, $b+\mu v+\ell=0$ and $-\gamma v+\beta\ell\mu^{-1}=0$.
Solving yields
\[
v=\frac{|b|}{\mu(1+\bet^{-1}\gam)}>0.
\]
The main heuristic step is again to substitute $v$ in for $V(s)$ in \eqref{23}.
This yields the approximation
\[
\tilde X(t)=X_0+\tilde bt+\sig W(t)+\tilde L(t),
\]
where
\[
\tilde b=\frac{b}{1+\beta\gamma^{-1}}.
\]
Thus $\tilde X$, that approximates $X$, is merely a reflected BM, and therefore
$E\tilde X(\iy)=\sig^2/(2|\tilde b|)$. Hence we obtain an approximation
$\widetilde{SD}$ for $SD$ in the form
\begin{equation}\label{1b}
\widetilde{SD}=1+E\tilde X(\iy)=1+\frac{\sig^2}{2|\tilde b|}= 1+\frac{\sig^2}{2|b|}(1+\beta\gamma^{-1}).
\end{equation}

As in the case of \S \ref{sec41},
the approximation is expected to improve as the mean length of the vacations grows,
for the same reasons.

\subsection{Numerics}

In both cases,
a standard Euler-Maruyama method is used to simulate the SDE \cite{higham}.
The Skorokhod constraining mechanism,
associated with the boundary term $L$, is treated by projecting $X$ back to zero whenever its
iteration assumes a negative value. The birth-death processes are treated by drawing
Bernoulli RVs, with suitable state dependent bias, to dictate upward and downward jumps.

The time step parameter is denoted by $\del$, the number of steps by
$N$, and the length of the simulated time interval is thus given by $T=N\del$.

\begin{figure}[p]
\center{\footnotesize
\includegraphics[width=0.48\textwidth,height=0.4\textwidth]{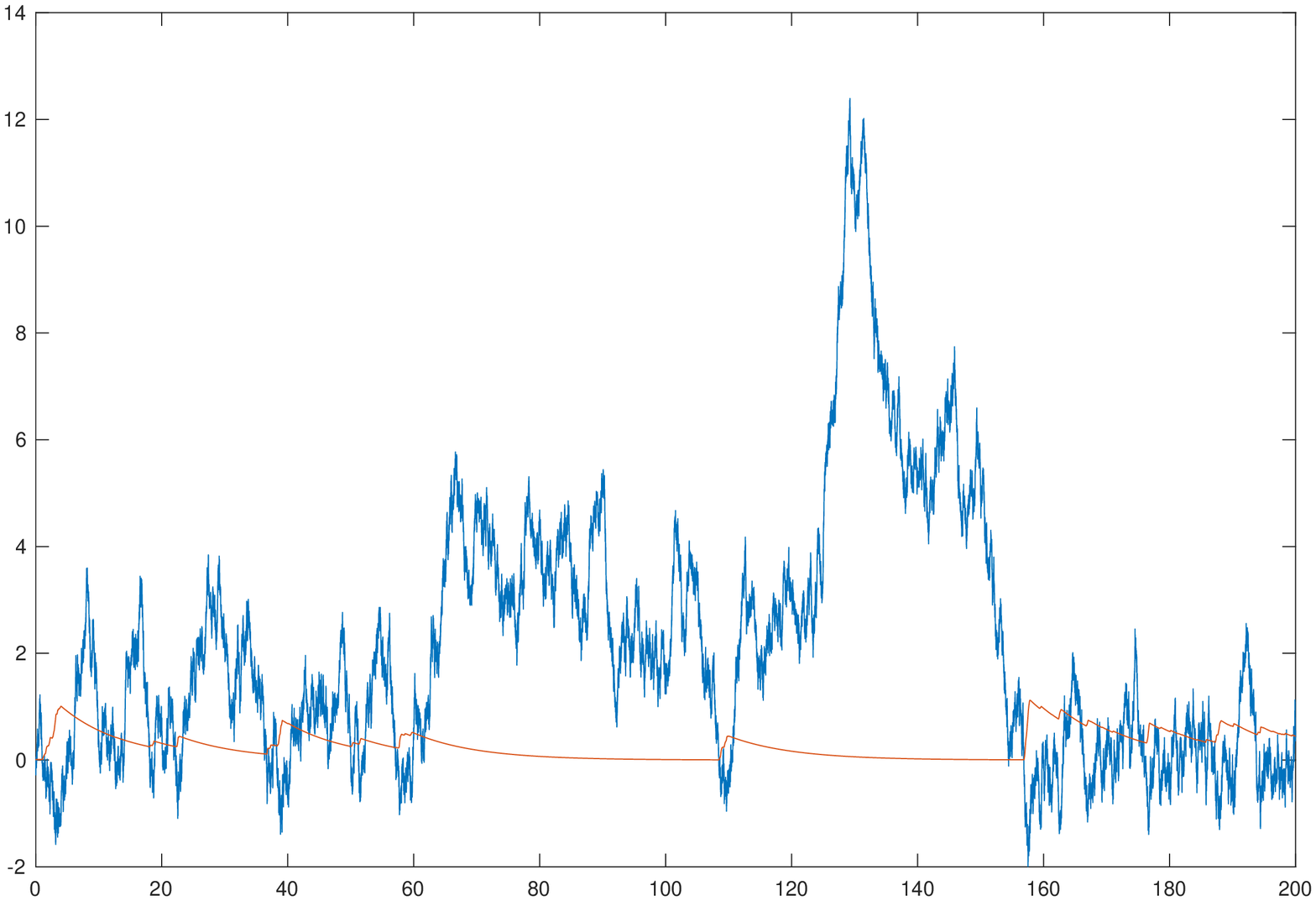}
\includegraphics[width=0.48\textwidth,height=0.4\textwidth]{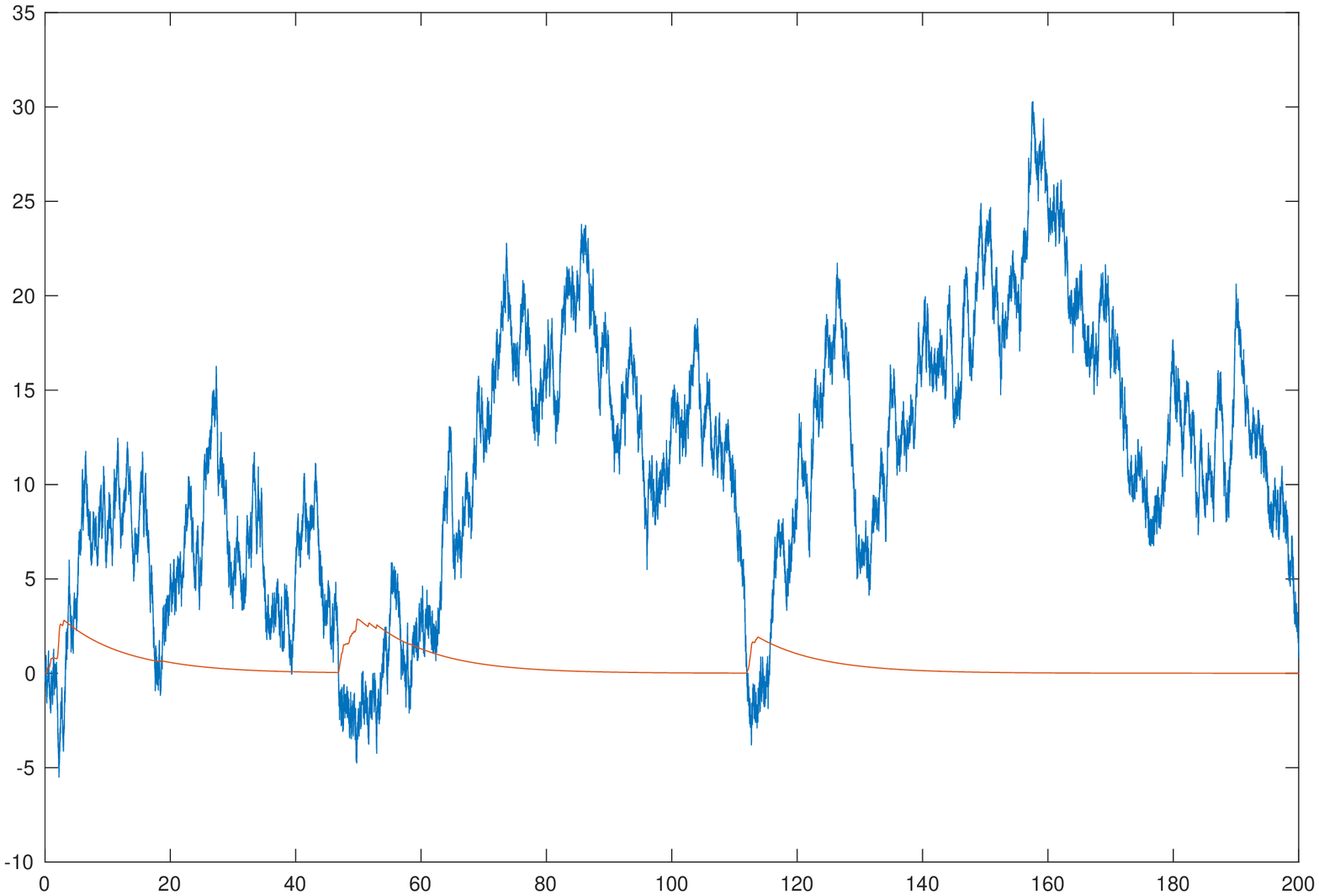}
\\[-1.5em]
(a) HW, $\sig=1$ \hspace{17em} (b) HW, $\sig=3$\\
\includegraphics[width=0.48\textwidth,height=0.4\textwidth]{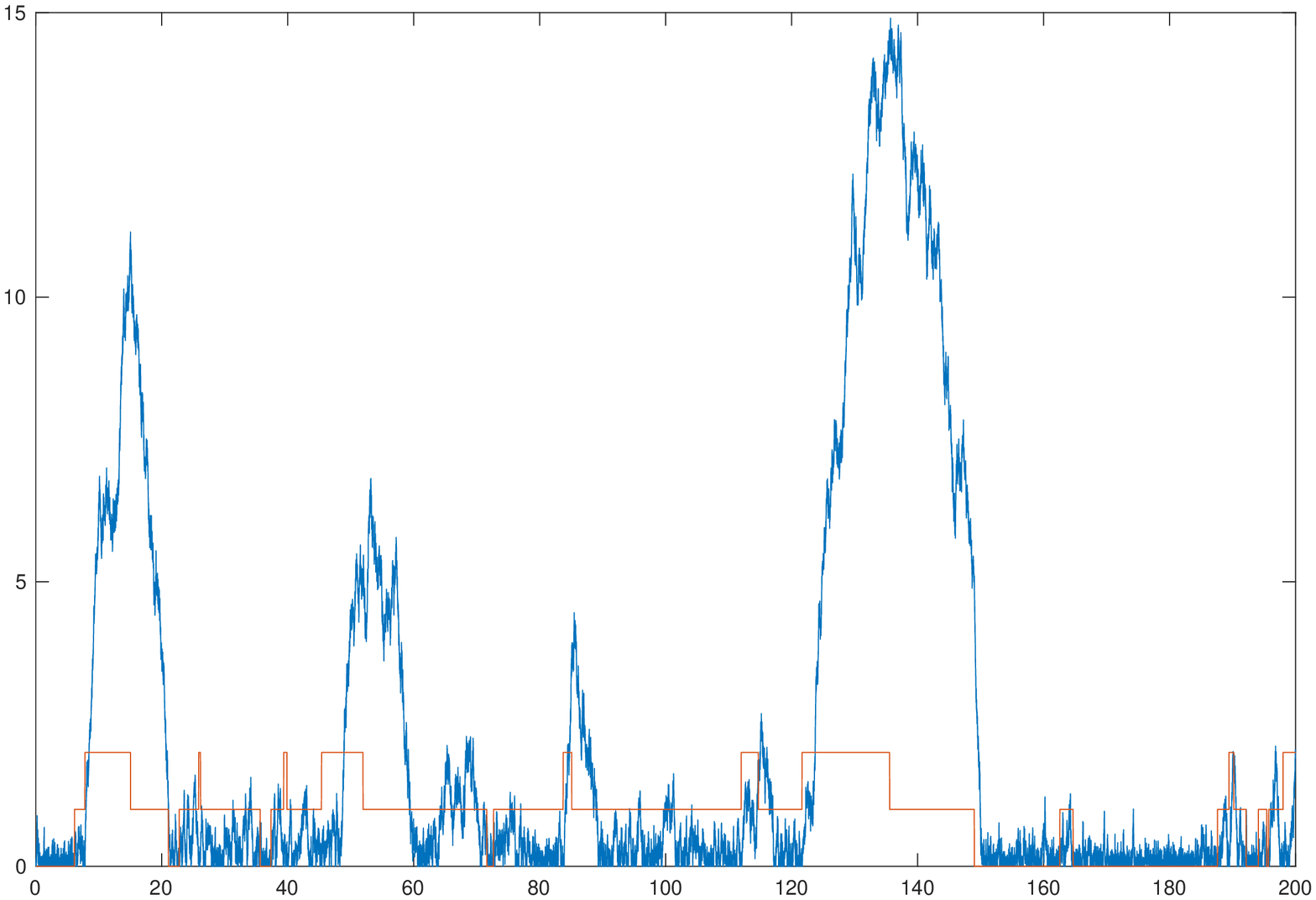}
\includegraphics[width=0.48\textwidth,height=0.4\textwidth]{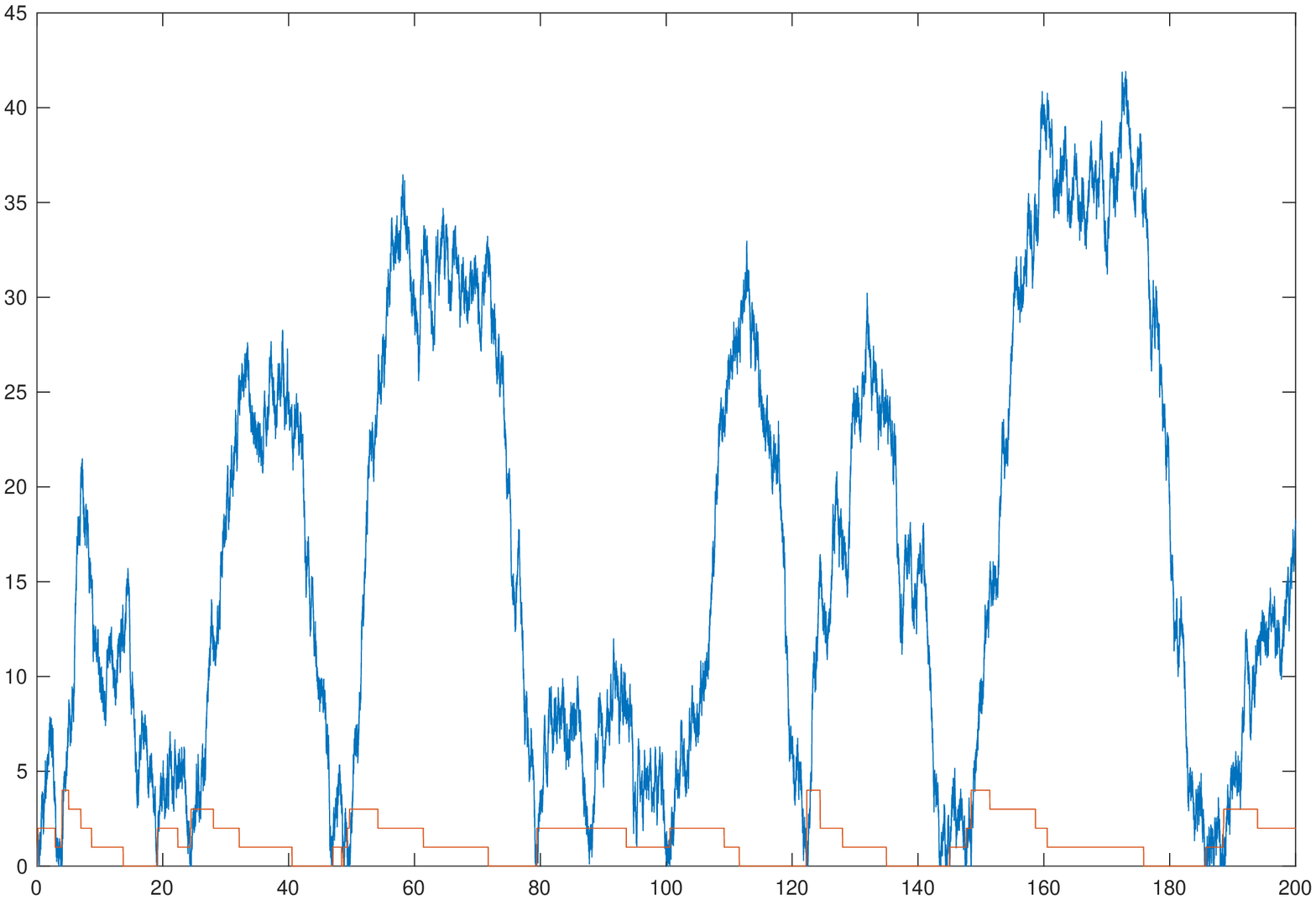}
\\[-1.5em]
(c) NDS, $\sig=1$ \hspace{16em} (d) NDS, $\sig=3$\\
}
\caption{\sl\footnotesize Parts (a) and (b) [resp., (c) and (d)] show sample paths
of the pair $X$ (blue) and $V$ (orange)
corresponding to the HW [resp., NDS] limit law, for different values of $\sig$.
}
\label{fig1}
\end{figure}

\paragraph{Sample paths}
First we present some sample paths of each of the coupled pairs \eqref{21} and \eqref{23},
where the coupling between the processes $X$ and $V$ is apparent.
The time step parameter is taken as $\del=10^{-3}$.
Here the number of steps is $N=2\times10^5$,
corresponding to a time interval of length $T=200$.
Sample paths for equation \eqref{21}, corresponding to the HW regime,
are shown in parts (a) and (b) of Figure \ref{fig1},
where $\sig=1$ and $\sig=3$, respectively (the remaining parameters are taken to be
$b=-0.3$, $\beta=2$, $\gamma=0.1$ $\mu=2$).
In both (a) and (b) we notice clearly the effect of $X$ on $V$.
Namely, an increase of $V$ occurs when $X+V$ is negative.
It is also noticeable that
$X$ reaches greater values in (b), where the diffusion coefficient is greater, than in (a).

For the NDS case, the sample paths of equation \eqref{23}
are shown in parts (c) and (d) of Figure \ref{fig1}.
Again, $\sig$ takes the two values $1$ and $3$, respectively
(and the remaining parameters are now $b=-3$, $\beta=2$, $\gamma=0.1$ $\mu=2$).
Here, the effect of each process on the other is visible.
Upward jumps of $V$ occur only when the diffusion process visits zero,
and its downward jumps occur only on excursions of $X$ away from zero.
The effect of $V$ on $X$ is particularly sharp in (c):
on each excursion of $X$, the path has strong tendency to increase when
$V=2$, but it decreases rapidly when $V=1$.
A similar dependence of the structure of the excursions on the value of $V$
occurs in (d), where $X$ has a strong negative drift when $V$ becomes zero.
Finally, as in the previous case, when $\sig$ is greater (that is, in (d)), $X$ reaches
higher values on excursions than when $\sig$ is smaller.

\paragraph{The HW regime}
We now use simulations to estimate the quality of the heuristic prediction
\eqref{1a} under various conditions.

As reference for the level of accuracy, we consider the model with
$V=0$ for which a theoretical value of $POW_0$ is known,
and compare to it simulation runs. This we do for the set of parameters
$\beta=2$, $\gamma=0.1$, $b=-2$, $\mu=1$, $\sigma=3$.
The theoretical value of $POW_0$ (given by formula \eqref{120}) and the
simulation results of $8$ runs with $N=2\times10^8$ steps (for the case with $V=0$,
i.e., based on sample paths of \eqref{114})
are summarized in Table \ref{table1}.
\begin{table}[h]
\center
\begin{tabular}{c || c | c | c | c | c | c | c | c || c}
$POW_0$ &
1 & 2 & 3 & 4 & 5 & 6 & 7 & 8 & max.\ dev.
\\ \hline
$0.2470$
&
$0.2487$ & $0.2473$ & $0.2487$ & $0.2441$ & $0.2437$ & $0.2459$ &
$0.2456$ & $0.2459$
&
$0.0033$
\end{tabular}
\caption{\sl\footnotesize
Simulation results for estimating $POW_0$ for $N=2\times10^8$ steps.
The theoretical value is shown on the left. The maximum absolute deviation
from the theoretical result is shown on the right.
}
\label{table1}
\end{table}

Among these 8 runs, the maximum absolute deviation
away from the theoretical value is  $0.0033$,
that is, less than $0.4\%$. This figure is sufficient for the purposes
of this study, and therefore in the simulations described below
we keep this value of $N$.

We control the length of vacations by varying $\gamma$.
The larger $\gamma$ is, the shorter is the
expected length of vacations.
Since the heuristic is based on substituting the long run average of $V$
for $V$ in the $X$ dynamics,
it is expected that for long vacations (small $\gamma$) the heuristic
provides accurate predictions. We also recall that when $\gamma=\mu$,
the equation for $X$ decouples from that of $V$ and the heuristic prediction
is exact.

\begin{figure}[p]
\center{\footnotesize
\includegraphics[width=0.48\textwidth,height=0.4\textwidth]{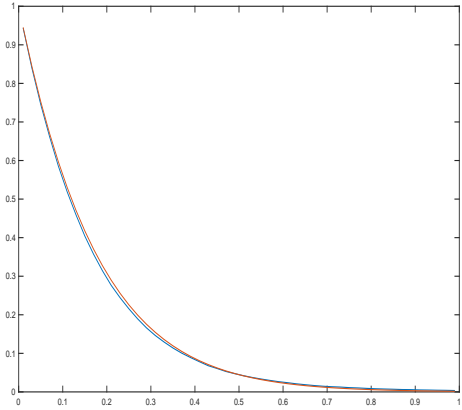}
\includegraphics[width=0.48\textwidth,height=0.4\textwidth]{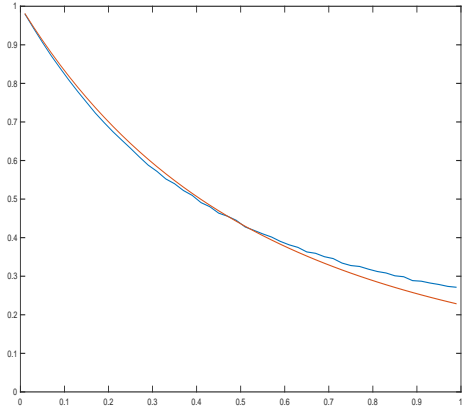}
\\[-1.5em]
(a) $\mu=0.5$, $\sig=1$ \hspace{17em} (b) $\mu=0.5$, $\sig=3$\\
\includegraphics[width=0.48\textwidth,height=0.4\textwidth]{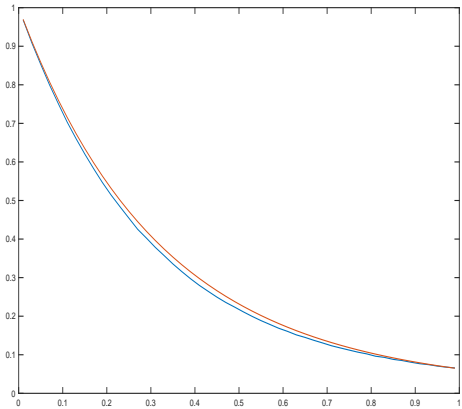}
\includegraphics[width=0.48\textwidth,height=0.4\textwidth]{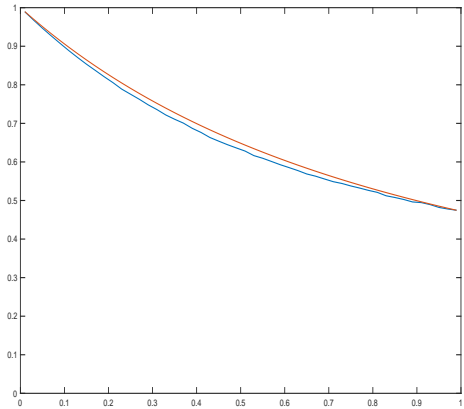}
\\[-1.5em]
(c) $\mu=1$, $\sig=1$ \hspace{16em} (d) $\mu=1$, $\sig=3$\\
}
\caption{\sl\footnotesize $POW$ obtained by simulation (blue),
and $\widetilde{POW}$ of formula \eqref{1a} (orange)
as a function of $\gamma$ in the range $0.01$ to $1$,
for different values of $\mu$ and $\sigma$.
Maximum absolute error: (a) $0.013$, (b) $0.042$, (c) $0.017$, (d) $0.020$.
}
\label{fig2}
\end{figure}

The results of our simulation are shown in Figure \ref{fig2}.
These four graphs show the simulation results of $POW$ and
the heuristic prediction $\widetilde{POW}$ of formula \eqref{1a}
as a function of $\gamma$, for two values of $\mu$
and two values of $\sigma$.
Specifically, $\gamma$ ranges between $0$ and $1$,
and the remaining parameters are taken as $b=-2$, $\beta=2$,
and $\mu\in\{0.5,1\}$ and $\sigma\in\{1,3\}$.

The arguments given above suggest that the graphs of $POW$ and $\widetilde{POW}$
should meet at two points, namely $\gamma=0$ and when $\gamma=\mu$.
This is seen very clearly in all parts of Figure \ref{fig2}.
As for the level of accuracy, it is overall very good in
cases (a), (c) and (d), and is somewhat less satisfactory in case (b).
For the actual numerical values of the maximal error sizes, see the figure description.
Overall, in all these cases the error is no greater than $5\%$.

Finally, the general behavior observed, where $POW$ (simulated and predicted)
is decreasing as $\gamma$ increases, is explained by the fact that
when the vacations are long ($\gamma$ small),
the system has effectively less service capacity,
consequently it is more loaded, and the probability of wait must increase.

\paragraph{The NDS regime}
In the case of the NDS regime, the simulations are aimed at testing the accuracy
of the prediction of formula \eqref{1b} for the slowdown.

Again we start by considering the reference model with $V=0$, for which there is a precise formula
for the slowdown. The parameters are taken to be
$\beta=5$, $\gam=3$, $b=6$, $\mu=2$, $\sig=3$.
It turns out that $\del$ must be calibrated.
With $\del=10^{-3}$ and $N=10^8$, there is a significant bias between the simulation and
theoretical value, explained by the inaccuracy introduced by
the constraining mechanism at zero, as an approximation for the boundary term $L$.
Whereas the size of this error converges to zero as $\del\to0$ by theoretical results,
the actual error for the above value of $\del$ is too large for our purposes.
When we reduce to $\del=10^{-4}$ and keep $N=10^8$, the bias is considerably smaller. The values
of 8 runs (based on simulating sample paths of \eqref{1c}) appear in Table \ref{table2},
along with the theoretical value (given by formula \eqref{1d})
and the maximal absolute error.
\begin{table}[h]
\center
\begin{tabular}{c || c | c | c | c | c | c | c | c || c}
$SD_0$ &
1 & 2 & 3 & 4 & 5 & 6 & 7 & 8 & max.\ dev.
\\ \hline
$1.7500$
&
$1.7288$ & $1.7340$ & $1.7352$ & $1.7220$ & $1.7410$ & $1.7231$ & $1.7382$ & $1.7402$
&
$0.0280$
\end{tabular}
\caption{\sl\footnotesize
Simulation results for estimating $SD_0$ for $N=10^8$ steps.
The theoretical value is shown on the left. The maximum absolute deviation
from the theoretical result is shown on the right.
}
\label{table2}
\end{table}

The maximal relative error is $0.0160$, that is less than $2\%$,
and is sufficient for our purposes.
In what follows we keep these values of $\del$ and $N$.

\begin{figure}[p]
\center{\footnotesize
\includegraphics[width=0.48\textwidth,height=0.4\textwidth]{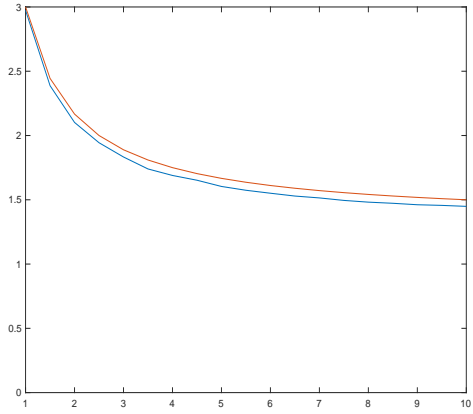}
\includegraphics[width=0.48\textwidth,height=0.4\textwidth]{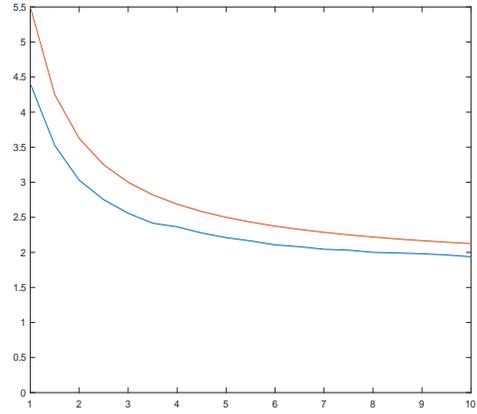}
\\[-1.5em]
(a) $\mu=2$, $\sig=2$ \hspace{17em} (b) $\mu=2$, $\sig=3$\\
\includegraphics[width=0.48\textwidth,height=0.4\textwidth]{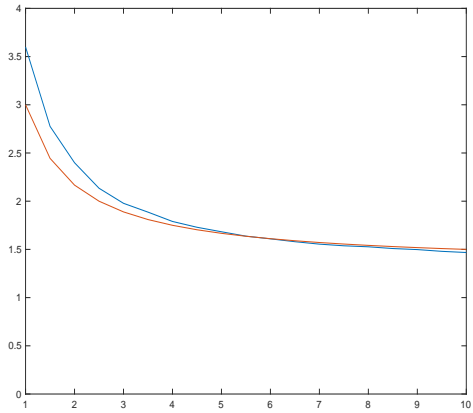}
\includegraphics[width=0.48\textwidth,height=0.4\textwidth]{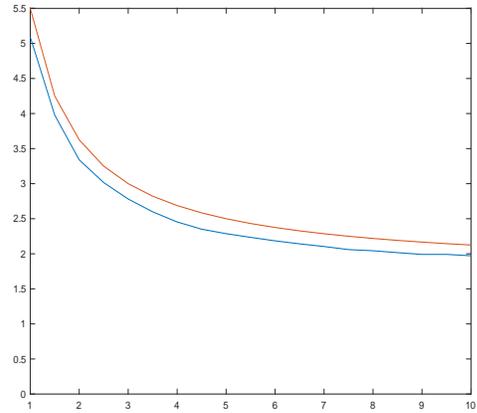}
\\[-1.5em]
(c) $\mu=4$, $\sig=2$ \hspace{16em} (d) $\mu=4$, $\sig=3$\\
}
\caption{\sl\footnotesize $SD$ obtained by simulation (blue),
and $\widetilde{SD}$ of formula \eqref{1a} (orange)
as a function of $\gamma$ in the range $1$ to $10$,
for different values of $\mu$ and $\sigma$.
Maximum absolute [resp., relative] error:
(a) $0.0696$ [$0.0390$],
(b) $1.0931$ [$0.1988$],
(c) $0.5996$ [$0.1999$],
(d) $0.4081$ [$0.0900$].
}
\label{fig3}
\end{figure}

Figure \ref{fig3} shows the simulation results of $SD$ and
the heuristic prediction $\widetilde{SD}$ of formula \eqref{1b}
as a function of $\gamma$, for two values of $\mu$
and two values of $\sigma$.
The parameters were chosen differently than in the HW case. Our concern
here was to calibrate the parameters so as to reach mean delay and mean service time of similar order.
This occurs when the slowdown is not too far from the value $2$.
Specifically, $\gamma$ ranges between $1$ and $10$,
and the remaining parameters are $b=-6$, $\beta=5$,
and $\mu\in\{2,4\}$ and $\sigma\in\{2,3\}$.

Overall, the accuracy of the heuristic prediction is worse than
in the HW case, with relative errors reaching as high as $20\%$
in some cases
(see description of Figure \ref{fig3}), although in parts of the ranges considered
the relative error is considerably smaller. In lack of better approximations,
these estimates may be useful in applications as first order approximations.

\appendix

\section{Appendix}\manualnames{A}
\subsection{One-dimensional Skorokhod problem}
\label{SP}

In this appendix we briefly review some well known properties of the one-dimensional Skorokhod problem.
For proofs of the results here, see \cite[Ch.\ 8]{ChungWilliams}.

\begin{definition}\label{def:sp}
	Given $x\in\D_+([0,\infty),\R)$ we say that a pair $(z,y)\in\D([0,\infty),\R_+)\times\D_0([0,\infty),\R_+)$ satisfies the one-dimensional Skorokhod problem for $x$ if the following conditions hold,
	\begin{itemize}
		\item[1.] $z(t)=x(t)+y(t)$ for all $t\geq0$;
		\item[2.] $y$ is non-decreasing and can only increase when $z$ is zero, i.e., $\int_0^t z(s)dy(s)=0$, $t\ge0$.
	\end{itemize}
\end{definition}

\begin{proposition}\label{prop:sp}
	Given $x\in\D_+([0,\infty),\R)$ there exists a unique solution $(z,y)$ of the one-dimensional Skorokhod problem for $h$ given by $(z,y)=(\Gam_1,\Gam_2)(h)$, where, for $t\geq0$,
	\begin{align}\label{eq:sm1}
	\Gam_1(x)(t)&=x(t)+\Gam_2(x)(t),\\ \label{eq:sm2}
	\Gam_2(x)(t)&=\sup_{0\leq s\leq t}(x(s))^-.
	\end{align}
	Consequently, the following properties hold:
	\begin{itemize}
		\item[1.] Oscillation inequality: given $x\in\D_+([0,\infty),\R)$ and $0\le s<t<\infty$,
		\begin{align}\label{eq:smosc}
		\text{Osc}(\Gam_1(x),[s,t])\leq\text{Osc}(x,[s,t])\qquad\text{and}\qquad\text{Osc}(\Gam_2(x),[s,t])\leq\text{Osc}(x,[s,t]).
		\end{align}
		\item[2.] Lipschitz continuity: for $x_1,x_2\in\D_+([0,\infty),\R)$ and $t\ge0$,
		\begin{align}
		\sup_{0\le s\le t}|\Gam_1(x_1)(s)-\Gam_1(x_2)(s)|&\leq2\sup_{0\le s\le t}|x_1(s)-x_2(s)|,\\
		\sup_{0\le s\le t}|\Gam_2(x_1)(s)-\Gam_2(x_2)(s)|&\leq\sup_{0\le s\le t}|x_1(s)-x_2(s)|.
		\end{align}
	\end{itemize}
\end{proposition}

\subsection{Nonexistence of relevant scaling for $\al\in[0,\frac{1}{2})$}
\label{sec:a2}

Here we provide an argument showing that for $\al$ in the range $[0,\frac{1}{2})$
there can be no rescaling of the server population process
under which the pair of processes (queue length, server population)
remains asymptotically coupled. This is argued by proving the following claim:
{\it Given any $T\in(0,\iy)$,
the {\it unnormalized} process $V^n$, if started at zero, remains zero on the interval
$[0,T]$ with probability tending to $1$ as $n\to0$.}

To prove the claim, let us first show that Lemma \ref{lem:eXn} remains valid for this range
of $\al$. By \eqref{eq:tildeX} and \eqref{16},
$0\le e^n_X(t)=n^{-1/2}(N^n-X^n(t))^+\le n^{-1/2}N^n\le n^{\al-1/2}\to0$.

Next consider the equation \eqref{11} for $V^n(t)$ with initial condition $V^n(0)=0$.
Let $s^n_1$ denote the first time when $V^n$ assumes the value $1$.
Our goal is to show that $P(s^n_1\le T)\to0$.

In equation \eqref{34}, the term $\int_0^tn^{-1}\mu^n\tilde N^n(s)ds$ vanishes
for all $t\le s^n_1$. The remaining terms in \eqref{34} are $C$-tight (recall
$e^n_X\To0$), and thus $\|\xi^n\|_{s^n_1\w T}$ is a tight sequence of RVs.
As a result of \eqref{30}, this is true also for $\|L^n\|_{s^n_1\w T}$.
By \eqref{28} and
\eqref{24} and the convergence $\mu^n/n$ to a positive constant,
we obtain that $k_n:=n^{-\al+\frac{1}{2}}\int_0^{s^n_1\w T}I^n(s)ds$ is
a tight sequence of RVs.
By the equation \eqref{11} for $V^n$, and the definition of $s^n_1$,
$S_B(\beta^n\int_0^{s^n_1}I^n(s)ds)=1$.
Consequently $\int_0^{s^n_1}I^n(s)ds\ge c\tau_1$, where $\tau_1$
is the first jump time of $S_B$, that is specifically an exponential with parameter 1.
Combining these facts,
\[
P(s^n_1\le T)\le P\lt(\int_0^{s^n_1\w T}I^n(s)ds\ge c\tau_1\rt)
=P(k_n\ge cn^{-\al+\frac{1}{2}}\tau_1)\to0,
\]
by the tightness of $k_n$.
\qed

\skp

\noi
{\bf Acknowledgment.}
The authors are grateful to Professor Ward Whitt for referring them to
\cite{whi-04} for formula \eqref{120}.
Research of RA supported in part by the ISF (grant 1184/16).
Research of DL supported in part by a Zuckerman fellowship.

\footnotesize

\bibliographystyle{is-abbrv}

\bibliography{bib}

\end{document}